\documentclass[12pt,draft,reqno]{amsart}
\usepackage{amssymb,delarray}
\usepackage{amsfonts}

\makeindex{}

\def\le{\leqslant}
\def\ge{\geqslant}

\newtheorem{thm}{Theorem}[section]
\newtheorem{lem}[thm]{Lemma}
\newtheorem{prop}[thm]{Proposition}
\newtheorem{claim}[thm]{Claim}

\newtheorem{cor}[thm]{Corollary}

\begin{document}
\baselineskip=14.pt plus 2pt 

\title[Old and new examples]{Old and  new examples of surfaces of general type
with $p_g=0$}
\author[Vik.S.~Kulikov]{Vik. S.~Kulikov}
\address{Steklov Mathematical Institute
} \email{kulikov@mi.ras.ru}

\dedicatory{} \subjclass{}
\thanks{The work  was partially supported by
RFBR ({\rm No.} 02-01-00786). 
}

\begin{abstract}
Surfaces of general type with geometric genus $p_g=0$, which can
be given as Galois covering of the projective plane branched over
an arrangement of lines with Galois group $G=(\mathbb Z/q\mathbb
Z)^k$, where $k\geq 2$ and $q$ is a prime number, are
investigated. The classical Godeaux surface, Campedelli surfaces,
Burniat surfaces, and a new surface $X$ with $K_X^2=6$ and
$(\mathbb Z/3\mathbb Z)^3\subset \mbox{Tors}\, (X)$ can be
obtained as such coverings. It is proved that the group of
automorphisms of a generic surface of the Campedelli type is
isomorphic to $(\mathbb Z/2\mathbb Z)^3$. The irreducible
components of the moduli space containing the Burniat surfaces are
described. It is shown that the Burniat surface $S$ with $K_S^2=2$
has the torsion group $\mbox{Tors}\, (S)\simeq (\mathbb Z/2\mathbb
Z)^3$, (therefore, it belongs to the family of the Campedelli
surfaces), i.e., the corresponding statement in the papers of C.
Peters "On certain examples of surfaces with $p_g=0$" in Nagoya
Math. J. {\bf 66} (1977), and I. Dolgachev "Algebraic surfaces
with $q=p_g=0$" in {\it Algebraic surfaces}, Liguori, Napoli
(1977), and in the book of W. Barth, C. Peters, A. Van de Ven
"Compact complex surfaces", p. 237, about the torsion group of the
Burniat surface $S$ with $K_S^2=2$ is not correct.
\end{abstract}

\maketitle
\setcounter{tocdepth}{2}


\def\st{{\sf st}}

\setcounter{section}{-1}
\section{Introduction}

As is known, the self-intersection number of the canonical class
of the surfaces of general type with geometric genus $p_g=0$ can
take the values $K^2=1,\dots, 9$, and in the past century the
existence of such surfaces for all possible values of $K^2$ was
proved. Nevertheless, our knowledge about the  surfaces of general
type with $p_g=0$ is far from completeness. In particular, the
moduli spaces of such surfaces are not described completely up to
now. Moreover, the list of all possible abelian groups, which can
be realized as the torsion group of such surfaces, is unknown.

In the paper we investigate surfaces of general type with $p_g=0$,
which can be given as Galois covering of the projective plane
branched over an arrangement of lines with Galois group
$G=(\mathbb Z/q\mathbb Z)^k$, where $k\geq 2$ and $q$ is a prime
number. In particular, the classical Godeaux surface \cite{God},
Campedelli surfaces \cite{Cam}, \cite{Mi}, Burniat surfaces
\cite{Bu}, and a new surface $X$ with $K_X^2=6$ and $$(\mathbb
Z/3\mathbb Z)^3\subset \mbox{Tors}\, (X)= \mbox{Tors}\,
H_1(X,\mathbb Z)  =\mbox{Tors}\, H^2(X,\mathbb Z)$$  can be
obtained as such coverings. It is proved that the group of
automorphisms of a generic surface of Campedelli type is
isomorphic to $(\mathbb Z/2\mathbb Z)^3$. It is shown that the
Burniat surface $S$ with $K_S^2=2$ has the torsion group
$\mbox{Tors}\, (S)\simeq (\mathbb Z/2\mathbb Z)^3$ (therefore, it
belongs to the family of the Campedelli surfaces (\cite{Mi}, see
also Proposition \ref{peters}), i.e., the corresponding statement
in \cite{Pet}, \cite{Dol}, and in \cite{B-P-V}, p. 237, about the
torsion group of the Burniat surface $S$ with $K_S^2=2$ is not
correct.

The irreducible components of  the moduli space containing the
Burniat surfaces are described. The description is depicted in the
following diagram

\begin{picture}(300,80)
\put(-5,60){$\mathcal M_2=\mathcal C$}
\put(3,20){$\widetilde{\mathcal B}_2$}\put(3,40){$\bigcup$}
\put(20,25){\vector(1,0){30}} \put(58,60){$\mathcal M_3$}
\put(60,20){$\widetilde{\mathcal B}_3$}\put(60,40){$\bigcup$}
\put(83,25){\vector(1,0){20}} \put(110,60){$\mathcal
M_4^{\prime\prime}\subset \mathcal M_4=\mathcal
M_4^{\prime\prime}\cup \mathcal M_4^{\prime}$}
\put(113,20){$\widetilde{\mathcal B}_4^{\prime\prime}\, \, \subset
\, \widetilde{\mathcal B}_4\, \, =\, \widetilde{\mathcal
B}_4^{\prime\prime}\, \, \sqcup \, \, \widetilde{\mathcal
B}_4^{\prime}$}\put(113,40){$\bigcup$}\put(150,40){$\bigcup$}
\put(238,25){\vector(1,0){20}} \put(270,60){$\mathcal M_5$}
\put(272,20){$\widetilde{\mathcal B}_5$}\put(272,40){$\bigcup$}
\put(293,25){\vector(1,0){20}} \put(320,60){$\mathcal M_6$}
\put(322,20){$\widetilde{\mathcal B}_6$}\put(322,40){$\bigcup$}

\end{picture}
\newline where  $\mathcal M_k$ is the union of irreducible components of
the moduli space of surface of general type with $p_g=0$ and
$K^2=k$ containing the Burniat surfaces and $\mathcal C$ is the
moduli space of the Campedelli surfaces. The points in the
subvarieties $\widetilde{\mathcal B}_k$ of $\mathcal M_k$
correspond to the Burniat surfaces. The varieties
$\widetilde{\mathcal B}_k$ for $k\neq 4$ are unirational  and
$\widetilde{\mathcal B}_4$ consists of two rational surfaces (the
points of the irreducible component $\widetilde{\mathcal
B}_4^{\prime\prime}$ parametrize the Burniat surfaces with $K^2=4$
having "$-2$"-curves), $\widetilde{\mathcal B}_2$ consists of a
single point, $\widetilde{\mathcal B}_3$ is a rational curve,
$\dim \widetilde{\mathcal B}_5=3$, and $\dim \widetilde{\mathcal
B}_6=4$. The subvarieties $\widetilde{\mathcal B}_k$ are
everywhere dense in $\mathcal M_k$ for $k\geq 4$, $\dim \mathcal
M_3=4$, and as is known (see \cite{Mi}),  $\mathcal C$ is
unirational, $\dim \mathcal C=6$. The arrows in the diagram show
the adjacency of the components (for example, $\widetilde{\mathcal
B}_3 \to \widetilde{\mathcal B}_4^{\prime\prime}$ means that the
Burniat surfaces with $K^2=3$ are degenerations of Burniat
surfaces with $K^2=4$ having "$-2$"-curves).

The paper is organized as follows. In section 1, we recall the
basic facts about Galois coverings $g:Y\to \mathbb P^2$ of the
plane $\mathbb P^2$ with Galois group $G=(\mathbb Z/q\mathbb Z)^k$
branched along a line arrangement $\overline L\subset \mathbb P^2$
and show how to obtain a resolution $X$ of the singular points of
$Y$ in terms of the singular points of $\overline L$. Then these
results are used in section 2 for the calculations of $K_X^2$ and
the topological Euler characteristic $e(X)$. In section 3, we
recall an algorithm of calculation of the geometric genus of $X$.
Section 4 is devoted to the examples and the result mentioned
above.

{\it Acknowledgement.} I would like to express my gratitude to the
University of Padova (Italy) for its hospitality during the early
stages of the preparation of this paper.

\section{Abelian coverings of the plane branched over an arrangement of lines}

By a Galois covering of a smooth algebraic variety $Y$ we mean a
finite morphism $h:X\to Y$ of a normal algebraic variety $X$ to
$Y$ such that the function fields imbedding $\mathbb C (Y)\subset
\mathbb C (X)$ induced by $h$ is a Galois extension. As is well
known, a finite morphism $h:X\to Y$ is a Galois covering with
Galois group $G$ if and only if $G$ coincides with the group of
covering transformations and the latter acts transitively on every
fiber of $h$. Besides, a finite branched covering is Galois if and
only if the un-ramified part of the covering (i.e., the
restriction to the complements of the ramification and branch
loci) is Galois. In addition, a branched covering is determined up
to isomorphism by its un-ramified part. Moreover, a covering map
from the unramified part of one branched covering to the
unramified part of another one induces a covering morphism between
these branched coverings if the extension of the morphism of
underlying varieties to the branch loci is given. Let us recall
also that an unramified covering is Galois with Galois group $G$
if and only if it is a covering associated with an epimorphism of
the fundamental group of the underlying variety to $G$, and, in
particular, the Galois coverings with abelian Galois group $G$ are
in one-to-one correspondence with epimorphisms to $G$ of the first
homology group with integral coefficients. All these results are
well known and their most nontrivial part can be deduced, for
example, from the Grauert-Remmert existence theorem \cite{G-R}.

In what follows we deal only with coverings of the complex
projective plane $\mathbb P ^2$ ramified over an arrangement of
lines $\overline L=L_1\cup \dots \cup L_n$. The simple loops
$\lambda _i, 1\le i\le n,$ around the lines $L_i$ generate
$H_1(\mathbb P ^2 \setminus \overline L ,\mathbb Z )\simeq \mathbb
Z ^{n-1}$. They are subject to the relation
$$\lambda _1+\dots +\lambda _n=0.$$

As for general abelian Galois coverings, a Galois covering $g:Y
\to \mathbb P ^2$ of $\mathbb P ^2$ with abelian Galois group $G$
branched along $\overline L$ is determined uniquely by an
epimorphism $\varphi : H_1(\mathbb P ^2 \setminus \overline L,
\mathbb Z )\to G$, and it exists for any such epimorphism. The
covering $g$ is branched along a line $L_i\subset \overline L$ if
and only if $\varphi (\lambda_i)\neq 0$ and, moreover, the
ramification index of $g$ along $L_i$ coincides with the order of
the element  $\varphi (\lambda_i)$ in $G$.

Since $H_1(\mathbb P ^2 \setminus \overline L,\mathbb Z )\simeq
\mathbb Z ^{n-1}$, there exists, in particular, an {\it universal
covering} $g_{u(m)}: Y_{u(m)}\to \mathbb P ^2$ corresponding to
the natural epimorphism $\overline \varphi : H_1(\mathbb P ^2
\setminus \overline L, \mathbb Z )\to H_1(\mathbb P ^2 \setminus
\overline L,\mathbb Z/m\mathbb Z )= H_1(\mathbb P ^2 \setminus
\overline L ,\mathbb Z )\otimes (\mathbb Z /m\mathbb Z )$. The
simplest example of such coverings is the following one.
\newline
{\bf Example}. Let $\overline L=L_0+L_1+L_2\subset \mathbb P^2$ be
given by equation $x_0x_1x_2=0$, where $(x_0:x_1:x_2)$ are
homogeneous coordinates of $\mathbb P^2$. It is easy to see that
the covering $g_{u(m)}:\mathbb P^2\to \mathbb P^2$ given by
$y_i=x_i^m$, $i=0,1,2$, is associated with the epimorphism
$$\overline \varphi : H_1(\mathbb P ^2 \setminus \overline L,
\mathbb Z )\simeq \mathbb Z^2\to (\mathbb Z/m\mathbb Z )^2.$$
\vspace{0.1cm}

The following statement is an immediate consequence of the general
results on branched coverings mentioned above.

\begin{prop} \label{0.1} If $g: Y \to \mathbb P ^2$ is a Galois
covering with Galois group $G\simeq (\mathbb Z /m\mathbb Z )^k$
branched along $\overline L$, then $k\le n-1$ and for any
epimorphism $H_1(\mathbb P^2\setminus\overline L)\to G$ there
exists a unique Galois covering $h: Y_{u(m)}\to Y$ inducing this
epimorphism and such that $g_{u(m)}=g\circ h$\,. \qed
\end{prop}

In what follows we deal with Galois coverings with Galois group
$G\simeq (\mathbb Z /q\mathbb Z )^k$, where $q$ is a prime number,
and we construct them in a way described in the above proposition.

Put
$$G_u=\{ \, \overline \gamma= (\gamma_1,\dots ,\gamma _{n-1}) \, \mid \,
\gamma_i\in  \mathbb Z/q\mathbb Z \}\simeq (\mathbb Z /q\mathbb Z
)^{n-1} $$ and let $\check{G_u}\simeq (\mathbb Z /q\mathbb Z
)^{n-1}$ be the dual (as a vector space over $\mathbb Z/q\mathbb
Z$) group, the pairing $(\overline \gamma, \overline a)$ is given
by
$$(\overline \gamma,\overline a)=\sum_{j=1}^{n-1}\gamma _ja _j\in \mathbb
Z/q\mathbb Z$$ for $\overline
\gamma=(\gamma_1,\dots,\gamma_{n-1})\in G_n$ and $\overline
a=(a_1,\dots,a_{n-1})\in \check{G_n}$.

Without loss of generality we can assume that the universal
covering $g_{u}: Y_{u}\to \mathbb P ^2$ is associated with the
epimorphism $\overline \varphi : H_1(\mathbb P ^2 \setminus
\overline L , \mathbb Z )\to G_u$ sending $\lambda _n$ to
$(q-1,\dots , q-1)$ and $\lambda _i$ with $1\leq i\leq n-1$ to
$(0,\dots,0,1,0,\dots,0)$ with $1$ in the $i$-th place. We choose
an additional line $L_{\infty}\subset \mathbb P ^2$ in general
position with respect to $\overline L$ and introduce affine
coordinates $(x,y)$ in $\mathbb C ^2=\mathbb P ^2\setminus
L_{\infty}$. Let $l_i(x,y)=0$ be a linear equation of $L_i \cap
\mathbb C ^2$. Put $z_i=(l_il_{n}^{q-1})^{1/q}$, $1\leq i \leq
n-1$. Then the function field $K_{u}=\mathbb C (Y_{u})=\mathbb C
(x,y,z_1,\dots ,z_{n-1})$ of a normal variety $Y_{u}$ is the
extension of the function field $K=\mathbb C (x,y)$ of $\mathbb P
^2$ of degree $q^{n-1}$. (In other words, the pull-back of
$\mathbb P^2\setminus L_\infty$ in $Y_{u}$ is naturally isomorphic
to the normalization of the affine subvariety of $\mathbb C^{n+1}$
given in coordinates $x,y,z_1,\dots, z_{n-1}$ by equations
$z_1^q=l_1l_n^{q-1},\dots , z_{n-1}^q=l_{n-1}l_n^{q-1}$.)

For a multi-index $\overline a=(a _1,\dots ,a_{n-1})$,
 $0\leq a _i\leq q-1$, we put
$$z^{\overline a}=\prod_{i=1}^{n-1}z_i^{a _i}.$$
The action of $\overline \gamma= (\gamma_1,\dots ,\gamma
_{n-1})\in G_u$ on $K_{u}$ is given by
$$\overline \gamma (z^{\overline a})=\mu ^{(\overline \gamma,\overline a)}
z^{\overline a},$$ where $\mu =e^{2\pi \sqrt{-1}/q}$ is the $q$-th
root of the unity. Therefore, we have $\mbox{Gal}(K_u/\mathbb
C[x,y])=G_u$ and
$$K_{u}=\bigoplus_{
0\le a _i\le q-1}\mathbb C (x,y)z^{\overline a}$$ is a
decomposition of the vector space $K_{u}$ over $\mathbb C (x,y)$
into a finite direct sum of degree $1$ representations of $G_u$.

Let $\varphi : H_1(\mathbb P^2 \setminus \overline L,\mathbb Z
)\to (\mathbb Z /q\mathbb Z)^k$ be an epimorphism given by
$\varphi (\lambda _i)=(a _{i,1},\dots ,a _{i,k})$, where $a
_{1,j}+\dots +a _{n,j}\equiv 0\, \text{mod}\, q$ for every
$j=1,\dots , k$, and let $g:Y \to \mathbb P ^2$ be the
corresponding Galois covering. The epimorphism $\varphi$ induces
the epimorphism $\psi :G_u\to G$. By Proposition \ref{0.1}, there
exists a unique Galois covering $h: Y_{u}\to Y$. It determines the
inclusion $h^{*} :\mathbb C (Y)\to K_{u}$ of the function field
$\mathbb C (Y)$ of $Y$ into the function field $K_{u}=\mathbb
C(Y_{u})$.

Since $\mbox{Gal}(K_{u}/h^*(\mathbb C (Y)))=\ker \psi$, obviously,
the field $h^*(\mathbb C(Y))$ coincides with the subfield
$K_{\varphi}=\mathbb C (x,y,w_1,\dots ,w_{k})$ of $K_{u}$, where
$w_j=z_1^{a _{1,j}}\cdot .\, .\, . \cdot z_{n-1}^{a_{n-1,j}}$, and
$$\mbox{Gal}(K_{u}/K_{\varphi})= \{ \,
(\gamma_1,\dots ,\gamma _{n-1})\in G\, \mid \, \sum_{i=1}^{n-1}
a_{i,j}\gamma_{i}\equiv 0\, 
(q),\, 1\leq j \leq k\, \}.$$

By construction, $Y$ is a normal surface with isolated
singularities. The singular points of $Y$ can appear only over the
$r$-fold points of $\overline L$ with $r\ge 2$, i.e., over
intersection points on $r$ lines $L_{i_1},\dots ,L_{i_r}$ of the
arrangement.

In what follows we call $2$ elements of $(\mathbb Z/q\mathbb Z)^k$
{\it linear independent over} $\mathbb Z/q\mathbb Z$ if they
generate in $(\mathbb Z/q\mathbb Z)^k$ a subgroup isomorphic to
$(\mathbb Z/q\mathbb Z)^2$.

\begin{lem} \label{lem 1.1} Let $p$ be a $2$-fold point of $\overline L$ and
$\varphi (\lambda_{i_1})$ and $\varphi (\lambda_{i_2})$ are linear
independent over $\mathbb Z /q\mathbb Z $ in $(\mathbb Z /q\mathbb
Z )^k$. Then the surface $Y$ is non-singular at each point of
$g^{-1}(p)$.
\end{lem}

\proof Let $p= L_{i_1}\cap L_{i_2}$. Choose a small neighborhood
$U$ of $p$ in $\mathbb P ^2$ and local analytic coordinates
$u_1,u_2$  in $U$ such that $U\simeq \{ \, |u_1|^2+|u_2|^2<
\varepsilon \, \}$ and $u_j=0$ is an equation of $L_{i_j}$. Then,
$H_1(U\setminus (L_{i_1}\cup L_{i_2}),\mathbb Z )\simeq \mathbb
Z\oplus \mathbb Z$. At any point $\widetilde p\in g^{-1}(p)$ the
germ $V\to U$ of the covering $Y\to\mathbb P^2$ is a
$G'$-covering, where $G'$ is the image of $H_1(U\setminus
(L_{i_1}\cup L_{i_2}),\mathbb Z)$ under the composition $\varphi
\circ i_*$ of $\varphi$ with the inclusion homomorphism $i_* :
H_1(U\setminus (L_{i_1}\cup L_{i_2}),\mathbb Z )\to H_1(\mathbb P
^2\setminus \overline L,\mathbb Z )$. Moreover, this $G'$-covering
is uniquely determined by $\varphi\circ i_*$. Identifying $\varphi
(\lambda _{i_1})$, $\varphi (\lambda _{i_2})$ with the standard
generators of $(\mathbb Z/q\mathbb Z)^2$ we get an isomorphism
between $V\to U$ and the covering determined by equations
$z_1^q=u_1, z_2^q=u_2$. Thus, $V$ is nonsingular. \qed

In our further examples, to resolve the singularities of $Y$ over
the $r$-fold points of $\overline L$ with $r\geq 3$, we blow up
all these points. Let $\sigma :\widetilde {\mathbb P ^2} \to
\mathbb P ^2$ be this blow up, $L'_i$ the strict transform of $L
_i$, $E_p$ the rational curve blown up over a $r$-fold point $p$,
and $\varepsilon_{p}\in H_1(\widetilde{\mathbb P ^2}\setminus
\sigma ^{-1}(\overline L),\mathbb Z ) =H_1(\mathbb P ^2 \setminus
\overline L,\mathbb Z )$ a simple loop around $E_{p}$.

The identification $H_1(\widetilde{\mathbb P ^2} \setminus \sigma
^{-1}(\overline L),\mathbb Z ) =H_1(\mathbb P ^2 \setminus
\overline L ,\mathbb Z )$ composed with $\varphi$ provides an
epimorphism $\varphi : H_1(\widetilde{\mathbb P ^2} \setminus
\sigma ^{-1}(\overline L),\mathbb Z )\to (\mathbb Z /q\mathbb Z
)^k$. Let consider the associated Galois covering $f:X \to
\widetilde{\mathbb P ^2}$.

The proof of the following statements is straightforward.

\begin{lem} \label{lem 1.2} Let $p=L_{i_1}\cap \dots \cap L_{i_r}$ be an
$r$-fold point of $\overline L$. Then $\varepsilon_{p}=\lambda
_{i_1}+\dots + \lambda _{i_r}$. \end{lem}

\proof    To establish the relation given by the Lemma, it is
sufficient to consider a generic line pencil passing through $p$.
\qed

\begin{lem} \label{lem 1.3} If for each $r$-fold point $p=L_{i_1}\cap \dots
\cap L_{i_r}$ of $\overline L$ with $r\geq 3$ either the pairs
$\varphi(\varepsilon_{p})$ and $\varphi(\lambda _{i_j})$, $1\leq j
\leq r$, are linear independent over $\mathbb Z /q\mathbb Z $ in
$(\mathbb Z /q\mathbb Z )^k$ or $\varphi(\varepsilon_{p})=0$, then
$X$ is nonsingular.
\end{lem}

\proof   It follows from Lemmas \ref{lem 1.1} and \ref{lem 1.2}.
\qed

Let $p_1,\dots, p_s$ be the set of $r$-fold points of a line
arrangement $\overline L$, $r\geq 2$, and let $\varphi :
H_1(\mathbb P^2 \setminus \overline L,\mathbb Z )\to (\mathbb Z
/q\mathbb Z)^k$ be an epimorphism given by $\varphi (\lambda
_i)=(a _{i,1},\dots ,a _{i,k})$, where $a _{1,j}+\dots +a
_{n,j}\equiv 0\, \text{mod}\, q$ for every $j=1,\dots , k$. Assume
that all singular points of $\overline L$ are  $\varphi$-{\it good
points}, i.e., for all $r$-fold points
$p_{i_1,\dots,i_r}=L_{i_1}\cap \dots \cap L_{i_r}$ of $\overline
L$ with $r\geq 2$,  either the pairs
$\varphi(\varepsilon_{p_{i_1,\dots,i_r}})$ and $\varphi(\lambda
_{i_j})$, $1\leq j \leq r$, are linear independent over $\mathbb Z
/q\mathbb Z $ in $(\mathbb Z /q\mathbb Z )^k$ or
$\varphi(\varepsilon_{p_{i_1,\dots,i_r}})=0$.  We say that a
$r$-fold point $p_{i_1,\dots,i_r}=L_{i_1}\cap \dots \cap L_{i_r}$
is {\it a non-branch point with respect to $\varphi$} if
$\varphi(\varepsilon_{p_{i_1,\dots,i_r}})=0$. Let $\sigma
:\widetilde {\mathbb P ^2} \to \mathbb P ^2$ be the blow up with
center at all $r$-fold points with $r\geq 3$  and at all $2$-fold
non-branch points of the arrangement $\overline L$. As a
consequence of Lemma \ref{lem 1.3}, the constructed surface $X$ is
a resolution of singularities of $Y$ and the covering $f$ is
included in the commutative diagram

\begin{picture}(280,90)
\put(130,70){$X$} 
\put(160,77){$\nu $} 
\put(143,75){\vector(1,0){40}} \put(192,70){$Y$} 
\put(133,65){\vector(0,-1){30}}\put(125,47){$f$} 
\put(195,65){\vector(0,-1){30}}\put(198,49){$g$} 
\put(130,20){$\widetilde{\mathbb P ^2}$} 
\put(160,17){$\sigma $} 
\put(145,25){\vector(1,0){40}} \put(191,20){$\mathbb P ^2$} 
\end{picture}
\newline inwhich $\nu$ is a regular birational map.

Let $N_{\varphi}$ be the set of all non-branch points with respect
to $\varphi$. Consider the subspace of $(\mathbb Z/q\mathbb
Z)^n=\{ (x_1,\dots,x_n) \, \mid \, x_i\in\mathbb Z/q\mathbb Z \}$
of solutions of the following system of linear  equations
\begin{equation} \label{equi}
\left\{
\begin{array}{l}
\displaystyle \hspace{0.5cm} \sum_{i=1}^n x_i=0,\vspace{0.3cm} \\
\displaystyle \sum_{i\in \{i_1,\dots,i_r\} }x_i=0
, \, \, \, \, \, p_{i_1,\dots,i_r}\in N_{\varphi}.
\end{array}
\right.
\end{equation}
Let
$n_{\varphi}$ be the rank of this linear system over $\mathbb
Z/q\mathbb Z $. We have $k\leq k_{\varphi}=n-n_{\varphi}$, since
the rank of the set of vectors $A_{\varphi}= \{
(a_{1,j}\dots,a_{n,j})\}_{\{j=1,\dots,k\}}$ is equal to $k$ and
the vectors from $A_{\varphi}$ satisfy equations (\ref{equi}). Let
us add $k_{\varphi}-k$ vectors to $A_{\varphi}$ to obtain  a basis
$A_{u,\varphi}$ over $\mathbb Z/q\mathbb Z$ of the space of
solutions of (\ref{equi})
$$A_{u,\varphi}=\{ (a_{1,j}\dots,a_{n,j})\}_{\{j=1,\dots,k_{\varphi}\}}$$
and consider the epimorphism $$\psi_{\varphi}: H_1(\mathbb P^2
\setminus \overline L,\mathbb Z )\to G_{u,\varphi}=(\mathbb Z
/q\mathbb Z)^{k_{\varphi}}$$ given by
$\psi_{\varphi}(\lambda_i)=(a _{i,1},\dots ,a _{i,k_{\varphi}})$.
Obviously, the epimorphism $\varphi$ can be decomposed into the
composition $\varphi=\eta\circ\psi_{\varphi}$, where $\eta
:(\mathbb Z/q\mathbb Z)^{k_{\varphi}}\to (\mathbb Z/q\mathbb Z)^k$
is the projection to the first $k$ coordinates. Let $\overline
f:\overline X\to \widetilde{\mathbb P}^2$ and
$h_{u,\varphi}:\overline X\to X$ be the Galois coverings
associated with $\psi_{\varphi}$ and $\eta$, respectively (see
Proposition \ref{0.1}). Note that the Galois group of the covering
$h_{u,\varphi}$ is isomorphic to $(\mathbb Z/q\mathbb
Z)^{k_{\varphi}-k}$.

The group $\mbox{Tors}(X)=\mbox{Tors}H_1(X,\mathbb Z)\simeq
\mbox{Tors}H^2(X,\mathbb Z)$ is called the {\it torsion group} of
$X$. Denote by $\mbox{Tors}_q(X)$ the subgroup of $\mbox{Tors}(X)$
consisting of the elements of order $q$.

The above consideration gives rise to the following
\begin{prop}\label{uni}
Let $f:X\to \widetilde{\mathbb P}^2$ be a Galois covering
associated with an epimorphism $\varphi : H_1(\mathbb P^2
\setminus \overline L,\mathbb Z )\to (\mathbb Z /q\mathbb Z)^k$
such that all singular points of the line arrangement $\overline
L$ are $\varphi$-good. Assume also that $\varphi (\lambda_i)\neq
0$ for each $L_i\subset \overline L$. Then
 $h_{u,\varphi}:\overline X\to X$ is unramified
covering.
\end{prop}

\begin{cor}\label{tors}
Let $f:X\to \widetilde{\mathbb P}^2$ be as in Proposition
\ref{uni}. If the irregularity $q(X)=\dim H^1(X,\mathcal O_X)=0$
and $k_{\varphi}-k>0$, then the $q$-torsion group $\mbox{\rm
Tors}_q(X)$ is non-trivial. In particular, there is an embedding
of $\ker \eta\simeq (\mathbb Z /q\mathbb Z)^{k_{\varphi}-k}$ to
$\mbox{\rm Tors}_q(X)$.
\end{cor}

\section{Calculation of $K^2$ and the Euler characteristic}

As above, let a Galois covering $g:Y \to \mathbb P ^2$ with Galois
group $G\simeq (\mathbb Z /q\mathbb Z)^{k}$ branched along a line
arrangement $\overline L=L_1+\dots +L_n$ be determined by an
epimorphism $\varphi : H_1(\mathbb P ^2 \setminus \overline L,
\mathbb Z )\to G$ such that $\varphi (\lambda_i)\neq 0$ for each
$L_i\subset \overline L$. Assume also that all singular points of
$\overline L$ are $\varphi$-good. Denote by $\sigma
:\widetilde{\mathbb P}^2\to \mathbb P^2$ the composition of
blowups with centers at the all $r$-fold points of $\overline L$
with $r\geq 3$ and at the all double points which are non-branch
points with respect to $\varphi$, and let $f:X\to
\widetilde{\mathbb P}^2$ be the covering induced by $\varphi$.
Since all singular points of $\overline L$ are $\varphi$-good, the
surface $X$ is non-singular.

Denote by $E_p=\sigma^{-1}(p)$ the curve blown up over a $r$-fold
point $p$, $L^{\prime}_i=\sigma^{-1}(L_i)$ the strict transform of
$L _i$, $C_i=f^{-1}(L^{\prime}_i)$, and $D_p=f^{-1}(E_p)$ the
strict transforms of the curves $L^{\prime}_i$ and $E_p$,
respectively. Let $T_r$ be the set of all $r$-fold points of
$\overline L$. Put
$$ T_r^{\prime}=\{ p\in T_r\, \, \mid \, \, p\, \, \mbox{is
a non-branch point of}\, \, \varphi \,\},$$
$T_r^{\prime\prime}=T_r \setminus T_r^{\prime}$, $\displaystyle
T^{\prime}=\bigcup_{r\geq 2}T_r^{\prime}$, $\displaystyle
T^{\prime\prime}=\bigcup_{r\geq 3}T_r^{\prime\prime}$, and
$T=T^{\prime}\cup T^{\prime\prime}$. Denote by $t_r^{\prime}=\#
T_r^{\prime}$ (respectively, $t_r^{\prime\prime}=\#
T_r^{\prime\prime}$) the number of the points belonging to
$T_r^{\prime}$ (respectively, $T_r^{\prime\prime}$) and put
$t_r=t_r^{\prime}+t_r^{\prime\prime}$. Note that the total
transform $f^*(L^{\prime}_i)=qC_i$ for each line $L_i\subset
\overline L$ and $f^*(E_p)=qD_p$ for each $p\in T^{\prime\prime}$.

\begin{thm} \label{K2}
The self-intersection number $K_X^2$ of the canonical class $K_X$
of $X$ is equal to
\begin{equation} \label{KK2}
\displaystyle  K_X^2=
q^{k-2}[(qn-n-3q)^2 -\sum_{r\geq 2}
(rq-q-r)^2t_r^{\prime}-\sum_{r\geq 3}
 (rq-2q-r+1)^2t_r^{\prime\prime} ].
\end{equation}

\end{thm}

\proof The canonical class of $\widetilde{\mathbb P}^2$ is equal
to $\displaystyle K_{\widetilde{\mathbb P}^2}=-3L+\sum_{p\in
T}E_p$, where $L=\sigma^*(\mathbb P^1)$ is the total transform of
a line $\mathbb P^1\subset \mathbb P^2$, and by adjunction
formula,
$$K_X=f^*(K_{\widetilde{\mathbb P}^2})+ (q-1)(\sum C_i+ \sum_{p\in
T^{\prime\prime}}D_p).$$ We have
$$\displaystyle q\sum C_i=f^*(nL
- \sum_{r\geq 3}\sum_{p\in T_r}rE_p-2\sum_{p\in T_2^{\prime}}E_p)
$$ and
$$\displaystyle q\sum_{p\in
T^{\prime\prime}}D_p=f^*(\sum_{p\in T^{\prime\prime}}E_p).$$
Therefore
$$\begin{array}{l}
\displaystyle qK_X=  qf^*(K_{\widetilde{\mathbb P}^2})+(q-1)(q\sum
C_i+
q\sum_{p\in T^{\prime\prime}}D_p)= \\
qf^*(-3L+\sum_{p\in T}E_p)+ \\ (q-1)f^*(nL - \sum_{r\geq
3}\sum_{p\in T_r}rE_p-2\sum_{p\in T_2^{\prime}}E_p)+ \\
(q-1)f^*(\sum_{p\in T^{\prime\prime}}E_p)
\end{array}
$$
and, finally,
$$ \begin{array}{l}   qK_X=  f^*((qn-n-3q)L- \\ \displaystyle
\mbox{\hspace{1.3cm}}\sum_{r\geq 2}\sum_{p\in
T_r^{\prime}}(rq-q-r)E_p-\sum_{r\geq 3}\sum_{p\in
 T_r^{\prime\prime}}(rq-2q-r+1)E_p).
\end{array}
$$ 
For each divisor $D\in \mbox{Pic}\, \widetilde{\mathbb P}^2$, we
have $$(f^*(D),f^*(D))_X=\deg f\cdot (D,D)_{\widetilde{\mathbb
P}^2}=q^k\cdot (D,D)_{\widetilde{\mathbb P}^2},$$ and Theorem
\ref{K2} follows from the equalities: $(L,L)_{\widetilde{\mathbb
P}^2}=1$, $(L,E_p)_{\widetilde{\mathbb P}^2}=0$ and
$(E_p,E_p)_{\widetilde{\mathbb P}^2}=-1$ for each $E_p$. \qed

In Section 4, we will apply Theorem \ref{K2} for the line
arrangements $\overline L$ and epimorphisms $\varphi$ with the
following properties: $t_2^{\prime}=0$, $t_4^{\prime}=0$, and
$t_r=0$ for $r\geq 5$. In this case formula (\ref{KK2}) takes the
following form
\begin{equation} \label{K2ex}
\displaystyle  K_X^2= q^{k-2}[(qn-n-3q)^2 - (2q-3)^2t_3^{\prime}-
 (q-2)^2t_3^{\prime\prime}- (2q-3)^2t_4^{\prime\prime}].
\end{equation}

Denote by
$$D_K=(qn-n-3q)L-\sum_{r\geq 2}\sum_{p\in
T_r^{\prime}}(rq-q-r)E_p-\sum_{r\geq 3}\sum_{p\in
 T_r^{\prime\prime}}(rq-2q-r+1)E_p.$$

Since $f$ is a finite Galois covering, we have the following
claim.
 \begin{claim} \label{ample} Let the divisor $D_K$ be big, i.e.,  $D_K^2>0$. Then
 \begin{itemize}
 \item[(i)] the surface $X$ is not minimal if and only if
 there is an irreducible curve $C\subset \widetilde{\mathbb P}^2$ such that
 $(D_K,C)_{\widetilde{\mathbb P}^2}<0$;
 \item[(i)] the canonical class of $X$ is not ample if and only if
 there is  an irreducible curve $C\subset \widetilde{\mathbb P}^2$ such that
 $(D_K,C)_{\widetilde{\mathbb P}^2}\leq 0$.
 \end{itemize}
 \end{claim}

\begin{thm}\label{euler}
The topological Euler characteristic of the surface $X$ is equal
to
\begin{equation} \label{eX}
\begin{array}{ll} \displaystyle
e(X)=& q^{k-2}(3q^2-2n(q^2-q)+q^2\sum_{r\geq 2}
t_r^{\prime}+(q-1)^2t_2^{\prime\prime} + \\
& ((r-1)(q-1)^2+1)\sum_{r\geq 3}t_r^{\prime\prime}).
\end{array}
\end{equation}
\end{thm}

\proof Denote by
$$B=\sum_{i=1}^nL_i^{\prime} +\sum_{p\in T^{\prime\prime}}E_p$$
the branch locus of $f$. It is easy to see that
\begin{equation} \label{??}
e(\mbox{Sing}B)=\# \mbox{Sing}B=t_2^{\prime\prime}+\sum_{r\geq
3}rt_r^{\prime\prime},
\end{equation}
where $\mbox{Sing}B$ is the
set of double points of the curve $B$.

The topological Euler characteristic of the curve $B$  is equal to
\begin{equation} \label{???}
e(B)= 2(n+\sum_{r\geq 3}t_r^{\prime\prime})-\# \mbox{Sing}B= 2n-
t_2^{\prime\prime}-\sum_{r\geq 3}(r-2)t_r^{\prime\prime},
\end{equation}
since $B$ is a divisor with normal crossings and the topological
Euler characteristic of each irreducible component of $B$ is equal
to $2$.

The topological Euler characteristic of $\widetilde{\mathbb P}^2$
is equal to
\begin{equation} \label{????}
\displaystyle e(\widetilde{\mathbb P}^2)=3+\sum_{r\geq 2}
t_r^{\prime}+\sum_{r\geq 3}
 t_r^{\prime\prime}.
\end{equation}

We have
\begin{equation} \label{?.?}
\begin{array}{ll}
\displaystyle e(X)= &
q^ke(\widetilde{\mathbb P}^2\setminus
B)+q^{k-1}e(B\setminus \mbox{Sing}B)+q^{k-2}e(\mbox{Sing}B)= \\
 & q^{k-2}(q^2e(\widetilde{\mathbb
 P}^2)-(q^2-q)e(B)-(q-1)e(\mbox{Sing}B)).
\end{array}
 \end{equation}
To complete the proof it is sufficient to substitute (\ref{??}) --
(\ref{????}) in (\ref{?.?}). \qed

For the line arrangements $\overline L$ and epimorphisms $\varphi$
with the following properties: $t_2^{\prime}=0$, $t_4^{\prime}=0$,
and $t_r=0$ for $r\geq 5$, formula (\ref{eX}) takes the following
form
\begin{equation} \label{eXex}
\begin{array}{ll} 
e(X)=  & q^{k-2}(3q^2-2n(q^2-q)+q^2
t_3^{\prime}+(q-1)^2t_2^{\prime\prime} +
\\ &
(2(q-1)^2+1)t_3^{\prime\prime}+(3(q-1)^2+1)t_4^{\prime\prime}).
\end{array}
\end{equation}

\section{Geometric genus calculation}
The aim of this section is to explain a general algorithm we will
use for  calculation of the geometric genus. In fact, if we
calculate the geometric genus of a covering, when we can calculate
its irregularity, since their difference is a topological
invariant equal to $\frac{K^2_X+e(X)}{12} -1$, due to Noether's
formula. In the calculation we use permanently the invariance of
the geometric genus under birational transformations, which allows
us at each step to use that nonsingular birational model which is
more convenient for the calculation.

The algorithm for calculation which we will use is by no means
new. It is contained, for example, in \cite{Kh-Ku}. Recall its
main steps.

{\bf 3.1.} {\it Reduction to cyclic coverings.} Let $g:Y_G\to
\mathbb P^2$, where $Y_G$ is supposed to be a normal surface, be a
Galois covering with abelian Galois group $G=(\mathbb Z /q\mathbb
Z)^k$ branched along curves $B_1,\dots , B_n\subset \mathbb P^2 $.
As above, such a covering is determined by an epimorphism $\varphi
:H_1(\mathbb P^2 \setminus \cup B_i) \to G$. Write it in a form
$$\varphi (\gamma _i)=
m_{1,i}\alpha_1+\dots +m_{k,i}\alpha_k, \qquad i=1,\dots ,\, n,$$
where $\alpha_j$ are standard generators of $G=\oplus(\mathbb
Z/q\mathbb Z)$, $\gamma _i$ are standard generators of
$H_1(\mathbb P^2 \setminus \cup B_i)$ dual to $B_i$ and
$m_{i,j}\in \mathbb Z /q\mathbb Z$, $0\leq m_{i,j} < q$, are
coordinates of $\varphi (\gamma _{i})$ with respect to $\alpha_j$.
In this notation, $Y_G$ is the normalization of the projective
closure of the affine surface $Y_{G,0}\subset \mathbb C^{m+2}$
given by
$$ z_j^{q}=\prod_{i=1}^n h_i^{m_{j,i}}(x,y), \qquad  \, j=1,\dots ,\, k,
$$
where $h_i(x,y)$ are  equations of $B_i$ in some chart $\mathbb
C^2 \subset \mathbb P^2$.

Let $X_G$ be the minimal desingularization of $Y_G$. As is known,
it exists, it is unique and the action of $G$ on $Y_G$ lifts
uniquely to a regular action on $X_G$.

Consider the action of $G$ on the space $H^0(X_G, \Omega
^2_{X_G})$ of regular 2-forms. It provides a decomposition
$$H^0(X_G,
\Omega ^2_{X_G}) =\oplus H_{(s_1,,\dots ,s_k)}$$ into the direct
sum of eigen-spaces $H_{(s_1,\dots, s_k)}$, where $\omega \in
H_{(s_1,\dots, s_k)}$ if and only if $\alpha_j(\omega )=e^{2\pi
s_j \sqrt{-1}/q}\cdot \omega $ for all $j=1,\dots , \, k$. Let
$H\subset G$ be a subgroup  and $G_1=G/H$. We have the following
commutative diagram

\begin{picture}(300,90)
\put(125,70){$Y_{G}$} 
\put(160,77){$h $} 
\put(143,75){\vector(1,0){40}} \put(192,70){$Y_{G_1}$} 
\put(131,65){\vector(0,-1){30}}\put(123,47){$g$} 
\put(195,65){\vector(0,-1){30}}\put(200,47){$g_1$} 
\put(128,20){${\mathbb P ^2}$} 
\put(160,13){$\mbox{id}$} 
\put(145,25){\vector(1,0){40}} \put(191,20){$\mathbb P ^2$} 
\end{picture}
\newline
where $f_1: Y_{G_1}\to \mathbb P^2$ is the Galois covering
corresponding to $\varphi _1=i\circ \varphi $ with $i:G\to
G_1=G/H$ being the canonical epimorphism. The map $h$ induces a
rational dominant (i.e., whose image is everywhere dense) map
$X_{G}\to X_{G_1}$, and the latter, as any rational dominant map
between nonsingular varieties, transforms holomorphic $2$-forms to
holomorphic $2$-forms. Thus, the subspace $h^*(H^0(
X_{G_1},\Omega^2_{X_{G_1}})) \subset H^0(X_G,\Omega^2_{X_G})$ is
well defined, and it coincides with the subspace
$$H^0(X_G,\Omega^2_{X_G})^H \subset
H^0(X_G,\Omega^2_{X_G})$$ of the elements fixed under the action
of $H$. On the other hand, an eigen-space $H_{(s_1,\dots, s_k)}$
is fixed under the action of $x_1\alpha_1+\dots x_k\alpha_k$ if
and only if $x_1s_1+\dots+x_ks_k=0\, (\mod  q)$. Hence, the sum
$\oplus H_{(\theta s_1,\dots, \theta s_k)}$ taken over
$\theta\in\mathbb Z/q\mathbb Z$ coincides with $H^0(\widetilde
X_G,\Omega^2_{\widetilde X_G})^H$, where
$$H=\{ \, x_1\alpha_1+\dots +x_k\alpha_k\, \, \vert \, \,
x_1s_1+\dots x_ks_k=0\, (q)\} .$$ So, this sum is isomorphic to
$H^0(X_{G/H},\Omega^2_{X_{G/H}}).$ These considerations give rise
to the following result.

\begin{prop} The geometric genus $p_g(X_G)=\dim
H^0(X_G,\Omega^2_{X_G})$ of $X_G$ is equal to
$$ p_g(X_G)= \sum_{H} p_g(X_{G/H}),$$
where the sum is taken over all subgroups $H$ of $G$ of $rk\,
H=rk\, G-1.$
\end{prop}

{\bf 3.2.} {\it Cyclic coverings.} Now, let $G=\mathbb Z /q\mathbb
Z $ be a cyclic group. To compute $p_g(X_G)$, let us choose
homogeneous coordinates $(x_0:x_1:x_2)$ in $\mathbb P^2$ such that
the line $x_0 =0 $ does not belong to the branch locus of $g:
Y_G\to \mathbb P^2$. As above, $Y_G$ is the normalization of the
projective closure of the surface in $\mathbb C^3$ given by
equation
$$z^p=h(x,y),$$
where $x=\frac{x_1}{x_0}$, $y=\frac{x_2}{x_0}$,
$$h(x,y)=\prod_{i=1}^n
h_i^{m_i}(x,y),$$ $h_i(x,y)$ are equations in $\mathbb C^2\subset
\mathbb P^2$ of the irreducible curves $B_i$ constituting the
branch locus, and $0<m_i<q$. Note that the degree
$$
\deg h(x,y)=\sum m_i\deg h_i(x,y)=mq
$$
is divisible by $q$, since the line $x_0 =0 $ does not belong to
the branch locus.

It is easy to see that over the chart $x_1\neq 0$ the variety
$Y_G$ coincides with the normalization of the surface in $\mathbb
C ^3$ given by equation
$$w^q=\widetilde h(u,v),$$
where $u=\frac{1}{x}$, $v=\frac{y}{x}$, $\widetilde
h(u,v)=u^{mq}h( \frac{1}{u},\frac{v}{u})$, and $w=zu^m$.

{\bf 3.2.1.} {\it Regularity condition over a generic point of the
base.} Consider $$\omega \in H^0(Y_G\setminus \text{Sing} Y_G,
\Omega ^2_{Y_G\setminus \text{Sing}Y_G})$$ and find a criterion of
its regularity outside the ramification and singular loci.

Over the chart $x_0\neq 0 $ the form $\omega$ can be written as
\begin{equation} \label{4} \omega =(\sum_{j=0}^{q-1}z^jg_j(x,y))\frac{dx\wedge
dy}{z^{q-1}},
\end{equation}
where $g_j(x,y)$ are rational
functions in $x$ and $y$. The form
$$\frac{dx\wedge dy}{z^{q-1}}$$
has neither poles nor zeros outside of the preimage of the branch
locus. Therefore, $\omega$ is regular at a point $(a,b)\not\in
\sum B_i$ if and only if all $g_j(x,y)$ are regular at the point.

In fact, if some $g_j(x,y)$ is not regular at $(a,b)$, then the
sum
$$\sum_{j=0}^{q-1}z^jg_j(x,y)$$ can be written
as
$$\frac{\sum_{j=0}^{q-1}z^jP_j(x,y)}{P_q(x,y)},$$
where $P_j(x,y)$, $j=0,\dots ,q$, are polynomials such that
$P_j(a,b)\neq 0$ for some $j<q$ and $P_q(a,b)=0$. Therefore,
$$\sum_{j=0}^{q-1}z^jP_j(a,b)=0$$
at all $q$ points belonging to $f^{-1}(a,b)$, since otherwise
$\omega$ would not be regular over $(a,b)$. On the other hand, it
is impossible, since a non-trivial polynomial of degree less than
$q$ can not have $q$ roots.

{\bf 3.2.2.} {\it Regularity condition over the line at infinity.}
Consider the form $\omega$ over the chart $x_1\neq 0$,
$$\omega = -(\sum_{j=0}^{q-1} w^j\frac{\widetilde g_j(u,v)}{u^{jm
+\deg g_j}}) \frac{1}{u^{3-m(q-1)}} \frac{du\wedge dv}{w^{q-1}}.$$
The similar arguments as above show that the regularity criterion
is equivalent to the following bound on the degrees of the
rational functions $g_j$
\begin{equation} \label{5} \deg g_j(x,y)\leq
(q-j-1)m-3.\end{equation}

{\bf 3.2.3}. {\it Regularity conditions over a nonsingular point
of the branch curve.} Consider our form
$$\omega =(\sum_{j=0}^{q-1}z^jg_j(x,y))\frac{dx\wedge dy}{z^{q-1}}$$
over a nonsingular point $(a,b)$ of one of the components,
$B_{i_0}$, of the branch curve. Let $r_j$ be the order of zero (or
of the pole if $r_j<0$) of the function $g_j$ along the curve
$B_{i_0}$, i.e., $g_j=\overline g_j\cdot h_{i_0}^{r_j}$ with
$\overline g_j$ having neither poles nor zeros along $B_{i_0}$.
Since $(a,b)$ is a nonsingular point of $B$, we can assume that
$h_{i_0}(x,y)$ and some function $g(x,y)$ are local analytic
coordinates  in some neighborhood $U$ of $(a,b)$ (denote them by
$u$ and $v$). So, over $U$ the surface $Y_G$ (after analytic
change of variables) is isomorphic to the normalization
$Y_{G,loc}$ of the surface in $\mathbb C ^3$ given by
$$z^{q}=u^{k_{i_0}}.$$
There is an analytic function $w$ in $Y_{G,loc}$ such that
$u=w^{q}$ and $z=w^{m_{i_0}}$, and such that $w$ and $y$ are
analytic coordinates in $Y_{G,loc}$. The differential 2-form
$\omega$ considered above has the following form in the new
coordinates
$$\omega =(\sum_{j=0}^{q-1}w^{jm_{i_0}}\overline g_j(x,y)w^{qr_j})
\frac{qw^{q-1}dw\wedge dv}{w^{(q-1)m_{i_0}}}.$$ It is easy to see
that
$$j_1m_{i_0}+qr_{j_1}+q-1-(q-1)m_{i_0}\neq
j_2m_{i_0}+qr_{j_2}+q-1-(q-1)m_{i_0}$$ if $0<m_{i_0}<q$, $0\leq
j_1, j_2\leq q-1$, and $j_1\neq j_2$. Therefore, $\omega$ is a
regular form  over a nonsingular point $(a,b)$ of $B_{i_0}$ if and
only if
$$jm_{i_0}+qr_{j}+q-1-(q-1)m_{i_0}\geq 0$$
for $0\leq j\leq q-1$. Moreover, if $\omega$ is a regular form
over $B_{i_0}$, then $r_j$ must be equal or greater than 0, since
for $0<m_{i_0}<q$, $0\leq j\leq q-1$, and $r_j\leq -1$, we obtain
that
$$jm_{i_0}+qr_{j}+q-1-(q-1)m_{i_0}< 0.$$
It follows that if $\omega$ is a regular form, then all rational
functions $g_j(x,y)$ are regular functions everywhere in $\mathbb
C ^2$ outside codimension 2. Thus, $g_j(x,y)$ should be
polynomials in $x$ and $y$. Moreover, the polynomials $g_j(x,y)$
must be divisible by $h^{r_j}_{i}(x,y)$, where $r_j$ is the
smallest non-negative integer satisfying the inequality
\begin{equation} \label{6}
qr_j\geq (q-j-1)m_i -q+1. \end{equation}

{\bf 3.2.4.} {\it Regularity conditions over singular points of
the branch curve.}  Let $\nu : X_G\to Y_G$ be the minimal
resolution of singularities of $Y_G$ and $E$ be the exceptional
divisor of $\nu$. Pick a composition $\sigma : \widetilde{\mathbb
P}^2\to \mathbb P^2$ of $\sigma$-processes with centers at
singular points of $B$ (and their preimages) such that $\sigma
^{-1}\circ f\circ \nu (E_i)$ is a curve for each irreducible
component $E_i$ of $E$. Let $Z$ be the normalization of
$\widetilde{\mathbb P}^2\times _{\mathbb P^2} Y_G$. Denote by $g:
X_G\to Z$ the bi-rational map induced by $\nu$ and $\sigma$. It
follows from the above choice of $\sigma$ that for any $\omega \in
H^0(Z\setminus \text{Sing}Z, \Omega ^2_{Z\setminus \text{Sing}Z})$
its pull-back $g^*(\omega )$ is regular at generic points of $E_i$
and, thus, extends to a regular form on the whole $X_G$. Hence,
$H^0(X_G, \Omega ^2_{X_G})$ is isomorphic to $H^0(Z\setminus
\text{Sing}\, Z, \Omega ^2_{Z\setminus \text{Sing}\, Z})$.

Therefore, it remains to consider a 2-form $\omega $ written as in
(\ref{4}) and to find a criterion of its regularity on
$Z\setminus\text{Sing}\, Z$. It can be done by performing, step by
step, the $\sigma$-processes chosen above. Let us accomplish only
the first step, since it is sufficient for the calculation in our
particular example which follows.

Represent, once more, $Y_G$ as normalization of the surface given
by
$$z^q=h(x,y).$$
Denote by $r$ the order of zero of $h(x,y)$ at the point $(0,0)$,
$r=sq+c$, $0\leq c <q$, and perform  the $\sigma$-process with
center at this point. In a suitable chart, this $\sigma$-process
$\sigma :\mathbb C ^2_{(u,v)}\to \mathbb C ^2_{(x,y)}$ is given by
$x=u, \, \, y=uv$. The normalization $Z_1$ of $Y_G\times _{\mathbb
C ^2_{(x,y)}} \mathbb C ^2_{(u,v)}$ is bi-rational to the
normalization of the surface given by
$$w^q=u^c\overline h(u,v),$$
where $w=z/u^s$ and $\overline h(u,v)=h(u,uv)/u^r$. We have
$$\begin{array}{ll}
\omega  & \displaystyle =(\sum_{j=0}^{q-1}z^jg_j(x,y))\frac{dx\wedge dy}{z^{q-1}}= \\
 & \displaystyle  =(\sum_{j=0}^{q-1}w^j\overline g_j(u,v)u^{sj+s_j+1-s(q-1)})
\frac{du\wedge dv}{w^{q-1}}, \end{array}
$$
where $s_j$ is the order of zero of $g_j(x,y)$ at $(0,0)$.
Applying (\ref{6}), we get necessary conditions for the regularity
of the pull-back of $\omega$ at generic points of the exceptional
divisor: the order of zero $s_j$ of each $g_j(x,y)$ at singular
point of the branch locus $B$ of order $r$ is the smallest integer
satisfying the inequality
\begin{equation} \label{7}
qs_j\geq
(q-j-1)r -2q+1. \end{equation}

To calculate the geometric genus of each $X_{G_i}$ we should find
explicitly all the regular $2$-forms, which are written as in
(\ref{4}) and satisfy criteria (\ref{5}) -- (\ref{7}).

The above discussion gives rise to the following statements for
$q=2$ and $3$ in the case of $q$-sheeted cyclic covering branched
along an arrangement of lines $\overline L=L_1+\dots+L_n$.
\begin{claim} \label{2cov}
Let $X$ be a desingularization of a double covering $g:Y\to\mathbb
P^2$ with branch locus $\overline L=L_1+\dots+L_n$. Denote by
$T_r$ the set of $r$-fold points of $\overline L$ and $T=\cup
T_r$. Then $n=2m$ is an even number and
$$
\begin{array}{l}
p_g(X)=  \dim_{\mathbb C} \{ \overline s\in H^0(\mathbb
P^2,\mathcal O_{\mathbb P^2}(m-3))\,  \mid \,\overline s \, \,
\mbox{has zero of}\,\, \mbox{order}\, \geq \\
\hspace{2.8cm} \lceil \frac{r+1}{2}\rceil -2 \, \, \mbox{at}\,\,
p\in T_r\, \mbox{for } \, \forall p\in T\} ,
\end{array}
$$
where $ \lceil \frac{a}{b} \rceil $ is the smallest  integer equal
or greater than $\frac{a}{b}$.
\end{claim}

\begin{claim} \label{3cov}
Let $X$ be a desingularization of a $3$-sheeted Galois covering
$g:Y\to\mathbb P^2$ in non-homogeneous coordinates given by
$$\displaystyle z^3=\prod_{i=1}^{n}l_i(x,y)^{m_i},$$
where $l_i(x,y)=0$ is an equation of $L_i$, each $m_i=1$ or $2$,
and $\sum m_i=3m$ is divisible by $3$. Denote by $T_r$ the set of
$r$-fold points of the divisor $\prod_{i=1}^{n}l_i(x,y)^{m_i}=0$,
$T=\cup T_r$, and $\widetilde
l(x,y)=\prod_{i=1}^{n}l_i(x,y)^{m_i-1}$. Then
$$
p_g(X)=  \dim_{\mathbb C}\mathcal P_0+ \dim_{\mathbb C}\mathcal
P_1 ,$$ where
$$
\begin{array}{l}
\mathcal P_0=\{  s\in H^0(\mathbb P^2,\mathcal O_{\mathbb
P^2}(2m-3-\sum (m_i-1)))\, \mid \, s\cdot\widetilde l \, \,
\mbox{has zero of}\,\, \\ \hspace{2.8cm} \mbox{order}\, \geq
2\lceil \frac{r+1}{3}\rceil -2 \, \, \mbox{at}\,\, p\in T_r\,
\mbox{for } \, \forall p\in T\}.
\end{array}
$$
and
$$
\begin{array}{l}
\mathcal P_1=\{ \overline s\in H^0(\mathbb P^2,\mathcal O_{\mathbb
P^2}(m-3))\, \mid \, \overline s \,\, \mbox{has zero of}\,\,
\mbox{order}\, \geq \lceil \frac{r+1}{3}\rceil -2
     \\
     \hspace{2.8cm}
\, \mbox{at}\,\, p\in T_r\, \mbox{for } \, \forall p\in T\}
\end{array}
$$
\end{claim}

\section{Examples}

{\bf 4.1.} {\it Campedelli surfaces}. Let $\overline L=L_1+ \dots
+L_7$ be a line arrangement in $\mathbb P^2$ consisting of seven
lines. We numerate them by the non-zero elements $\alpha_i\in
(\mathbb Z/2\mathbb Z)^3$. We will assume that the arrangement
$\overline L$ has not $r$-fold points with $r\geq 4$ and if
$\overline L$ has a tripe point $p_{\alpha_1,\alpha_2,\alpha_3}=
L_{\alpha_1}\cap L_{\alpha_2}\cap L_{\alpha_3}$, then
$\alpha_1+\alpha_2+\alpha_3\neq 0$. Such arrangement of lines is
called a {\it Campedelli arrangement}. Consider the covering $g:
Y\to \mathbb P^2$ induced by the epimorphism $\varphi :H_1(\mathbb
P^2\setminus \overline L, \mathbb Z)\to G=(\mathbb Z/2\mathbb
Z)^3$ given by $\varphi (\lambda_{\alpha_i})=\alpha_i$.

The surface $Y$ has the singular points lying only over the triple
points $p_{\alpha_1,\alpha_2,\alpha_3}$. To resolve them, let us
blow up the triple points and consider the induced Galois covering
$f: X \to \widetilde{\mathbb P}^2$, where $\sigma :
\widetilde{\mathbb P}^2\to {\mathbb P}^2$ is the composition of
blowups with centers at the triple points. We call the constructed
surface $X$ a {\it Campedelli surface}. Denote by
$E_{\alpha_i,\alpha_j,\alpha_k}=\sigma^{-1}(p_{\alpha_i,\alpha_j,\alpha_k})$
the exceptional curve lying over $p_{\alpha_i,\alpha_j,\alpha_k}$.
Since $\alpha_i+\alpha_j+\alpha_k\neq 0$ for triple points, each
curve $E_{\alpha_i,\alpha_j,\alpha_k}$ is a branch curve of the
covering $f$. It follows from Lemma \ref{lem 1.3} that $X$ is
non-singular, since $\varphi
(\varepsilon_{\alpha_i,\alpha_j,\alpha_k})=\alpha_i+\alpha_j+\alpha_k$
and $\alpha _i$ (respectively, $\alpha _j$ $\alpha _k$) are linear
independent in $G$. Indeed, $\alpha_i+\alpha_j+\alpha_k$ and
$\alpha _i$ are linear dependent if and only if
$\alpha_i+\alpha_j+\alpha_k=\alpha_i$, i.e., if and only if
$\alpha_j=\alpha_k$.

\begin{prop} \label{camp} The constructed Campedelli surfaces $X$ are
surfaces of general type with $K_X^2=2$, $p_g=0$, and
$\mbox{Tors}(X)=(\mathbb Z/2\mathbb Z)^3$.
\end{prop}

\proof Applying Claim \ref{ample}, we have $2K_X=\mid \widetilde
f^*(L)\mid $, where $L=\sigma^*(\mathbb P^1)$ is the total
transform of a line $\mathbb P^1\subset \mathbb P^2$. Therefore
$X$ is a surface of general type. Moreover, it is minimal, since
$(L,C)_{\widetilde{\mathbb P}^2}\geq 0$ for each curve $C\subset
\widetilde{\mathbb P}^2$. Applying (\ref{K2ex}) and (\ref{eXex}),
it is easy to see that $K_X^2=2$ and $e(X)=10$. Therefore, by
Noether's formula, $p_a=1-q+p_g=1$. As above, to calculate $p_g$,
it is enough to calculate the geometric genuses of  $7$ cyclic
coverings corresponding to $7$ epimorphisms $\psi_k$, $k=1,\dots
,7$, of $G=(\mathbb Z/2\mathbb Z)^3$ to the cyclic group $\mathbb
Z/2\mathbb Z$. It is easy to see that each of these coverings is
given  in non-homogeneous coordinates by the equation of the form
$w_k^2=l_{\alpha_{i_1}}l_{\alpha_{i_2}}l_{\alpha_{i_3}}l_{\alpha_{i_4}}$,
where $\alpha_{i_1},\alpha_{i_2},\alpha_{i_3},\alpha_{i_4}$ are
the elements of $G$ such that $\psi_k (\alpha_{i_j})=1$. Applying
Claim \ref{2cov}, one can easily check that the geometric genus of
each of these coverings is equal to zero. Thus, $X$ has the
geometric genus $p_g=0$.

To show that $\mbox{Tors}(X)=(\mathbb Z/2\mathbb Z)^3$, consider
the universal covering $f_{u(2)}: X_{u(2)}\to \widetilde{\mathbb
P}^2$ corresponding to the epimorphism
$$\overline
\varphi :H_1(\mathbb P^2\setminus \overline L, \mathbb Z)\to
H_1(\mathbb P^2\setminus \overline L, \mathbb Z/2\mathbb Z )\simeq
(\mathbb Z/2\mathbb Z)^6,$$ and the covering $h:X_{u(2)}\to X$
corresponding to an epimorphism $\psi :(\mathbb Z/2\mathbb Z)^6\to
G=(\mathbb Z/2\mathbb Z)^3$.
 By Proposition \ref{uni} and Corollary \ref{tors}, the covering $h$
is unramified and $(\mathbb Z/2\mathbb Z)^3\subset
\mbox{Tors}(X)$. Therefore, by  \cite{Mi},
$\mbox{Tors}(X)=(\mathbb Z/2\mathbb Z)^3$. \qed

The classical Campedelli surface $S$ (\cite{Cam}) is obtained as a
resolution of singularities of a double covering $\widetilde
g:\widetilde Y \to \mathbb P^2$ branched along the union of three
quadrics $Q_1,Q_2,Q_3$ and a quartic $C_4$ in $\mathbb P^2$ such
that the curve $B=Q_1+Q_2+Q_3+C_4$ has $6$ singular points of the
type $[3,3]$ (a singular point of the type $[3,3]$ means that
after the blow up with center at the singular point the strict
transform of the germ of $B$ consists of $3$ irreducible smooth
branches each pair of which meets transversally).

Let us show that the classical Campedelli surface $S$ is
isomorphic to $X$. Consider a covering $f:X\to \widetilde{\mathbb
P}^2$ branched along a Campedelli arrangement $\overline L=\sum
L_{\alpha_i}$, $\alpha_i\in (\mathbb Z/2\mathbb Z)^3\setminus \{
0\}$, having $3$ triple points. The arrangement $\overline L$ is
depicted in Fig. 1.

\begin{picture}(300,265)
\qbezier(100,200)(200,150)(300,100) \put(250,105){$L_{(0,0,1)}$}
\qbezier(50,40)(130,40)(270,40) \qbezier(55,30)(110,140)(165,250)
\qbezier(155,250)(210,140)(265,30)
\qbezier(160,260)(160,150)(160,30)
\qbezier(50,32)(130,99)(210,166) \qbezier(95,150)(160,107)(270,33)
\put(47,45){$p_2$} \put(265,45){$p_1$} \put(165,240){$p_3$}
\put(56,100){$L_{(1,0,0)}$} \put(162,70){$L_{(1,1,1)}$}
\put(192,190){$L_{(0,1,0)}$} \put(101,64){$L_{(0,1,1)}$}
\put(175,97){$L_{(1,0,1)}$} \put(101,29){$L_{(1,1,0)}$}
\end{picture}

\begin{center} Fig. 1 \end{center}
To see this isomorphism, let us consider the blowup $\sigma
:\widetilde{\mathbb P}^2 \to \mathbb P^2$ with center at the
points $p_1,p_2,p_3$ and denote by $E_i=\sigma^{-1}(p_i)$ the
exceptional curve lying over $p_i$, and the strict transforms
$\sigma^{-1}(L_{\alpha_i})\subset \widetilde{\mathbb P}^2$ we will
denote again by $L_{\alpha_i}$. One can check that
\begin{equation} \label{h}
\varphi (\varepsilon_i)=(0,0,1) \end{equation} for $i=1,2,3$.  The
curves $L_{(1,0,0)}$, $L_{(1,1,0)}$, $L_{(0,1,0)}$ in
$\widetilde{\mathbb P}^2$ have self-intersection numbers equal
$-1$. Therefore we can blow down them by monoidal transformation
$\tau :\widetilde{\mathbb P}^2\to \mathbb P^2$ (the composition
$\tau\circ\sigma^{-1} :\mathbb P^2\to \mathbb P^2$ is the
quadratic transformation of the plane with center at the points
$p_1,p_2,p_3$). The curves $L_i=\tau (E_i)$, $i=1,2,3$, and $\tau
(L_{(1,0,1)})$, $\tau (L_{(0,1,1)})$, $\tau (L_{(1,1,1)})$ are
lines and  $\tau (L_{(0,0,1)})$ is a conic in $\mathbb P^2$. We
have the following commutative diagram

\begin{picture}(300,90)
\put(130,70){$X$} 
\put(160,77){$\nu$} 
\put(143,75){\vector(1,0){40}} \put(192,70){$\check{Y}$} 
\put(133,65){\vector(0,-1){28}}\put(125,47){$f$} 
\put(195,65){\vector(0,-1){28}}\put(200,47){$\check{g}$} 
\put(130,20){$\widetilde{\mathbb P ^2}$} 
\put(160,15){$\tau$} 
\put(145,25){\vector(1,0){40}} \put(191,20){$\mathbb P ^2$} 
\end{picture}
\newline
where $\check{Y}$ is a normal surface, $\nu :X\to \check{Y}$ is a
bi-rational map and $\check{g}:\check{Y}\to \mathbb P^2$ is the
Galois covering branched along the curves $L_i$, $i=1,2,3$, $\tau
(L_{(1,0,1)})$, $\tau (L_{(0,1,1)})$, $\tau (L_{(1,1,1)})$, and
$\tau (L_{(0,0,1)})$. Since $\varphi (\lambda_{(0,0,1)})=(0,0,1)$,
and taking into account (\ref{h}) it is easy to see that
$\check{g}$ can be decomposed into the composition
$\check{g}=g_1\circ \widetilde g$, where $g_1:\mathbb P^2\to
\mathbb P^2$ is the Galois covering with Galois group
$G_1=(\mathbb Z/2\mathbb Z)^2$ branched along the lines $\tau
(L_{(1,0,1)})$, $\tau (L_{(0,1,1)})$, $\tau (L_{(1,1,1)})$ (see
Example in Section 1) and $\widetilde g:\widetilde Y\to \mathbb
P^2$ is the Galois covering with the Galois group $G_2=\mathbb
Z/2\mathbb Z$ branched along $Q_i=g_1^{-1}(L_i)$, $i=1,2,3$, and
$C_4=g_1^{-1}(\tau (L_{(0,0,1)}))$, where $Q_1,Q_2,Q_3$ are
quadrics and $C_4$ is a quartic in $\mathbb P^2$ such that the
curve $B=Q_1+Q_2+Q_3+C_4$ has $6$ singular points of the type
$[3,3]$. \qed

\begin{thm} \label{autcam} For a generic Campedelli surface $X$, the group
$\mbox{\rm Aut}\, (X)$ is isomorphic to $(\mathbb Z/2\mathbb Z)^3$
and coincides with the covering transformation group $G$ of
$f:X\to \widetilde{\mathbb P}^2$.
\end{thm}

\proof Let $\overline L\subset \mathbb P^2$ be a Campedelli
arrangement without triple points and assume that if an
automorphism $\widetilde h$ of $\mathbb P^2$ leaves fixed
$\overline L$ (i.e., $h(\overline L)=\overline L$), then
$h=\mbox{id}$. Consider the covering $g:Y\to \mathbb P^2$
associated with $\varphi :H_1(\mathbb P^2\setminus \overline L,
\mathbb Z)\to G=(\mathbb Z/2\mathbb Z)^3$ given by $\varphi
(\lambda_{\alpha_i})=\alpha_i$. Since the arrangement $\overline
L$ has not triple points, $Y=X$ is a nonsingular surface and
$g=f$. The morphism $f$ induces an extension of fields
$f^*(\mathbb C(\mathbb P^2))\subset \mathbb C(X)$. As in section
1, we choose a line at infinity $L_{\infty}$, coordinates $(x,y)$
in $\mathbb P^2\setminus L_{\infty}$ and identify $\mathbb
C(\mathbb P^2)$ with the field $\mathbb C(x,y)$ of rational
functions.

Consider an element $\alpha \in \mbox{\rm Tors}_2(X)=\mbox{\rm
Tors}\, (X)$, $\alpha\neq 0$. The linear system $|K_{X}+\alpha |$
is non-empty and $D\in |K_{X}+\alpha |$ for some $\alpha \in
\mbox{\rm Tors}_2(X_s)$ if and only if $2D= f^*(\widetilde L)$ for
some $\widetilde L \in | L|$, where $L$ is a line in $\mathbb
P^2$. Indeed, the linear system $|K_{X}+\alpha |$ is non-empty by
Riemann -- Roch Theorem, since $\dim H^2(X, \mathcal
O_{X}(K_{X}+\alpha))=\dim H^0(X, \mathcal O_{X}(\alpha))=0$. Let
$D_{\alpha}\in |K_{X}+\alpha |$. Then $2D_{\alpha}\in |2K_{X}|$.
By Riemann -- Roch Theorem, we have $ \dim
H^0(X,2K_{X})=K^2_{X}+1=3$. On the other hand,  it follows from
Claim \ref{ample} that $2K_{X}=f^*(L)$ and $ \dim H^0(\mathbb
P^2,\mathcal{O}_{\mathbb P^2}(1))=3$. Therefore, $|2K_{X}|
=f^*(|L|)$ and $D\in |K_{X}+\alpha |$ for some $\alpha \in
\mbox{\rm Tors}_2(X)$ if and only if $2D= f^*(\widetilde L)$ for
some $\widetilde L \in |L|$.

It is easy to see that there are exactly $7$ lines $\widetilde L
\in |L|$ for which the divisors $f^*(\widetilde L)$ are divisible
by $2$, namely, $L_{\alpha}\subset \overline L$, $\alpha \in
\mbox{\rm Tors}_2(X)$, $\alpha \neq 0$. So, we have
$\frac{1}{2}f^*(L_{\alpha})=D_{\alpha}\in |K_{X}+\alpha |$.

Let $h:X \to X$ be an isomorphism. Then it induces isomorphisms
$h^*: \mbox{\rm Tors}\, (X)\to \mbox{\rm Tors}\, (X)$ and
$$h^*: H^0(X, \mathcal
O_{X}(K_{X}+\alpha)) \to H^0(X, \mathcal
O_{X}(K_{X}+h^*(\alpha)))$$ for each $\alpha\in \mbox{\rm Tors}\,
(X)$. Therefore, $h^*(D_{\alpha})=D_{h^*(\alpha )}$ for $\alpha
\in \mbox{\rm Tors}_2(X)$, $\alpha \neq 0$. The automorphism $h$
induces the action $h^*$ on the group $\mbox{Div}\,X$ We have
$$
\begin{array}{l}
h^*(f^*(L_{\alpha_1}-L_{\alpha_2}))=h^*(2D_{\alpha_1}-2D_{\alpha_2})=2D_{h*(\alpha_1)}-2D_{h^*(\alpha_2)}=
\\ \hspace{3.8cm} f^*(L_{h^*(\alpha_1)}-L_{h^*(\alpha_2)})
\end{array}
$$
for any $\alpha_1,\alpha_2\in \mbox{\rm Tors}\, (X)$,
$\alpha_1\neq \alpha_2\neq 0$. Therefore,
\begin{equation}
\label{function}\displaystyle
h^*(f^*(\frac{l_{\alpha_1}(x,y)}{l_{\alpha_2}(x,y)}))=
c_{\alpha_1,\alpha_2}f^*(\frac{l_{h^*(\alpha_1)}(x,y)}{l_{h^*(\alpha_2)}(x,y)}),
\end{equation}
where $c_{\alpha_1,\alpha_2}$ is a constant, since any rational
function is defined uniquely up to multiplication by a constant by
its divisors of zeros and poles. It follows from (\ref{function})
that $h^*$ induces an automorphism $\widetilde h^*$ of $\mathbb
C(x,y)$ such that $f^*\circ\widetilde h^*=h^*\circ f^*$, since the
functions $\frac{l_{\alpha_1}(x,y)}{l_{\alpha_2}(x,y)}$ generate
the field $\mathbb C(x,y)$. Moreover, the automorphism $\widetilde
h^*$ induces an automorphism $\widetilde h$ of $\mathbb P^2$ such
that $\widetilde h(\overline L)=\overline L$. Therefore,
$\widetilde h=\mbox{\rm id}$ and $h\in \mbox{\rm Gal}\, (X/\mathbb
P^2)$. \qed
\begin{thm} \label{modcamp} {\rm ( cf. \cite{Mi} )}
The moduli space $\mathcal C$ of the Campedelli surfaces is an
unirational variety, $\dim \mathcal C=6$.
\end{thm}

\proof By the same arguments, that were used in the proof of
Theorem \ref{autcam}, one can show that two Campedelli surfaces
$X_1$ and $X_2$, defined by Campedelli line arrangements
$\overline L_1$ and $\overline L_2$, are isomorphic if and only if
there is a linear transformation $h$ of $\mathbb P^2$ sending
$\overline L_1$ to $\overline L_2$.

Applying a suitable linear transformation of $\mathbb P^2$ and a
suitable automorphism of $(\mathbb Z/2\mathbb Z)^3$,  we can
assume that for a line arrangement $\overline L=\sum L_{\alpha}$,
the lines $L_{(1,0,0)}$, $L_{(0,1,0)}$, $L_{(1,1,0)}$, and
$L_{(1,1,1)}$ are given respectively by $z_0=0$,  $z_1=0$,
$z_2=0$, and $z_0+z_1+z_2=0$. Therefore, a line arrangement
$\overline L$ is defined by a point in an everywhere dense subset
$V$ of $(\check{\mathbb P}^2\setminus \{ \mbox{ four points} \}
)^3$. Obviously, for any point $v_0\in V$, the set $A_{v_0}\subset
V$ consisting of the points for which the corresponding line
arrangements $\overline L_v$, $v\in A_{v_0}$, can be transformed
to $\overline L_{v_0}$ by linear transformations of $\mathbb P^2$,
is finite. Therefore, the moduli space $\mathcal C$ is an
unirational variety, $\dim \mathcal C=6$ (see also Corollaries
\ref{modBur} and \ref{cambur}).

 \qed \newline {\bf 4.2.} {\it Burniat
surfaces.} Let $\overline L_s=L_1+\dots +L_9$ be an arrangement in
$\mathbb P^2$ of nine lines depicted in Fig. 2. The arrangement
$\overline L_s$ has three $4$-fold points $p_1,p_2,p_3$ and $s$
($s=0,\dots, 4$) triple points $p_{3+i}$, $0< i\leq s$.

\begin{picture}(300,290)

\qbezier(50,40)(130,40)(270,40) \qbezier(55,30)(110,140)(165,250)
\qbezier(155,250)(210,140)(265,30)
\qbezier(160,280)(160,240)(160,200)
\qbezier(120,240)(160,240)(200,240)
\qbezier(155,250)(210,140)(265,30) \put(193,245){$L_9$}
\put(163,265){$L_8$} \qbezier(30,10)(50,30)(80,60)
\qbezier(20,66)(50,46)(80,26) \qbezier(240,53)(270,33)(300,13)
\qbezier(240,27)(270,47)(300,67) \put(53,47){$p_2$}
\put(260,48){$p_1$} \put(166,245){$p_3$} \put(99,150){$L_{7}$}
\put(40,10){$L_{5}$} \put(210,150){$L_{1}$} \put(280,10){$L_{3}$}
\put(285,47){$L_{2}$} \put(151,29){$L_{4}$} \put(27,47){$L_{6}$}
\end{picture}

\begin{center} Fig. 2 \end{center}
Such line arrangements we will call {\it Burinat arrangements}.
 Consider the covering $g: Y_s\to \mathbb P^2$ induced by the epimorphism
$\varphi :H_1(\mathbb P^2\setminus \overline L_s, \mathbb Z)\to
G=(\mathbb Z/2\mathbb Z)^2$ given by
\[
\begin{array}{l}
\varphi (\lambda_{1})=  \varphi (\lambda_{2})=\varphi
(\lambda_{3})=(1,0), \\
\varphi (\lambda_{4})=\varphi (\lambda_{5})=\varphi
(\lambda_{6})=(0,1), \\
\varphi (\lambda_{7})=\varphi (\lambda_{8})=\varphi
(\lambda_{9})=(1,1).
\end{array}
\]

The surface $Y_s$ has $3+s$ singular points lying over the
$4$-fold points $p_j$, $j=1,2,3$, and the triple points $p_{3+i}$,
$1\leq i\leq s$. Let $\sigma :\widetilde{\mathbb P}^2\to \mathbb
P^2$ be the composition of the blowups with centers at these
points. Denote by $E_j=\sigma^{-1}(p_j)$ the exceptional curve
lying over $p_j$, $1\leq j \leq 3+s$. Consider the induced Galois
covering $f: X_s \to \widetilde{\mathbb P}^2$. We have
\[
\begin{array}{c}
\varphi (\varepsilon_{1})= (1,1), \, \, \,  \varphi
(\varepsilon_{2})= (1,0), \, \, \, \varphi
(\varepsilon_{3})=(1,1), \\
\varphi (\varepsilon_{3+i})= (0,0)\, \, \mbox{for}\, \, 1\leq
i\leq s.
\end{array}
\]

Therefore, the curves $E_j$ are the branch curves of $f$ for
$j=1,2,3$ and the curves $E_{3+j}$ are not the branch curves of
$f$ for $j\geq 1$. Thus, by Lemma \ref{lem 1.3}, $X_s$ is a smooth
surface. Note that the number of triple and $4$-fold points of
$\overline L$ is less than $8$ and each $4$ of such points do not
lie in the same line. Therefore, $\widetilde{\mathbb P}^2$ is a
del Pezzo surface possibly  with "$-2$"-curves.

\begin{prop} \label{burnia}
The constructed above surfaces $X_s$ (they are called the Burniat
surfaces) are surfaces of general type with $K_{X_s}^2=(6-s)$ and
$p_g=0$.
\end{prop}

\proof By Claim \ref{ample}, we have $$2K_{X_s}=\mid \widetilde
f^*(3L-\sum_{i=1}^{3+s} E_i)\mid ,$$ where $L=\sigma^*(\mathbb
P^1)$ is the preimage of a line $\mathbb P^1\subset \mathbb P^2$.
Therefore, $X_s$ is a surface of general type and it is minimal,
since $\widetilde{\mathbb P}^2$ is a del Pezzo surface possibly
with "$-2$"-curves. Applying (\ref{K2ex}) and (\ref{eXex}), it is
easy to see that $K_{X_s}^2=6-s$ and $e(X)=6+s$. Therefore, by
Noether's formula, $p_a=1-q+p_g=1$. As above, to calculate $p_g$,
it is enough to calculate the geometric genera of the
desingularizations $\overline Z_i$ of $3$ cyclic coverings
$g_i:Z_i\to \mathbb P^2$ corresponding to $3$ epimorphisms
 from the group $G=(\mathbb Z/2\mathbb Z)^2$ to the cyclic group $\mathbb Z/2\mathbb Z$:

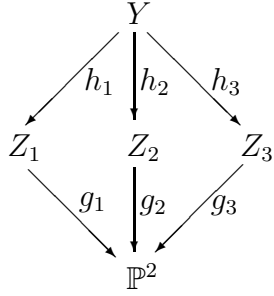
\begin{picture}(300,125)
\put(150,110){$Y$} \put(105,60){$Z_1$} \put(150,60){$Z_2$}
\put(193,60){$Z_3$}\put(150,10){$\mathbb P^2$}
\put(147,107){\vector(-1,-1){35}}\put(158,107){\vector(1,-1){35}}
\put(153,107){\vector(0,-1){33}}
\put(113,55){\vector(1,-1){33}}\put(194,57){\vector(-1,-1){33}}
\put(153,56){\vector(0,-1){32}}
\put(134,85){$h_1$}\put(155,85){$h_2$}\put(182,85){$h_3$}
\put(132,41){$g_1$}\put(155,40){$g_2$}\put(182,40){$g_3$}
\end{picture}
\newline where $g=g_i\circ h_i$ for $i=1,2,3$.
 These coverings are given in non-homogeneous coordinates respectively by the following
equations:
\begin{equation} \label{w1w2w3}
\begin{array}{l}
w_1^2=l_1l_2l_3l_4l_5l_6; \\ w_2^2=l_4l_5l_6l_7l_8l_9; \\
w_3^2=l_1l_2l_3l_7l_8l_9.
\end{array}
\end{equation}

Applying Claim  \ref{2cov}, one can easily check that the
geometric genus of each of these coverings vanishes, since each of
the arrangements given respectively by $l_1l_2l_3l_4l_5l_6=0$,
$l_4l_5l_6l_7l_8l_9=0$, and $l_1l_2l_3l_7l_8l_9=0$ has a $4$-fold
point. Thus, $X_s$ has the geometric genus $p_g=0$. \qed

Denote  by $\widetilde L_j$ the strict transform
$\sigma^{-1}(L_j)$ of the curve $L_j\subset \overline L_s$. Then
the divisor $\sum \widetilde L_j+\sum_{i=1}^3E_i$ is the branch
locus of the covering $f$. Put
\begin{equation} \label{notation}
\begin{array}{ll} 2C_j=f^{*}(\widetilde L_j), & j=1,\dots,9,
\\
2D_i=f^{*}(E_i) &  i=1,2,3,
\\ \hspace{0.1cm} D_k=f^{*}(E_k) &  3<k\leq 3+s,
\end{array}
\end{equation}
and denote by $t(L_j)$ the number of singular points of the
arrangement $\overline L_s$ lying on the line $L_j$.
\begin{claim} \label{k2} We have
\begin{itemize}
\item[(i)] the curves $C_j$, $j=1,\dots,9$, and $D_i$, $i=1,\dots,3+s$,  are non-singular;
\item[(ii)]  the geometric genus of a curve  $C_j$,
$j=1,\dots,9$, is equal to $g(D_i)=3-t(L_j)$;
\item[(iii)]  the geometric genus of a curve $D_i$
 is equal to $g(D_i)=1$ if $i=1,2,3$ and $g(D_i)=0$ if $3<i\leq 3+s$;
\item[(iv)] the self-intersection number of a curve  $C_j$,
$j=1,\dots,9$, is equal to $(C^2_j)_{X_s}=1-t(L_j)$;
\item[(v)] the self-intersection number of a curve $D_i$  is equal to
$(D^2_i)_{X_s}=-1$ if $i=1,2,3$ and $(D^2_i)_{X_s}=-4$ if $3<i\leq
3+s$.
\end{itemize}
\end{claim}
\proof Statement (i) is obvious.

(ii) -- (iii) The restriction of the covering map $f$ to a curve
$C_j$, $j=1,\dots,9$, is a two-sheeted covering of a rational
curve branched at $8-2t(L_j)$ points. Therefore,
$g(C_j)=3-t(L_j)$.

The restriction of $f$ to a curve $D_i$, $i=1,2,3$, is a
two-sheeted covering of a rational curve branched at $4$ points.
Therefore $g(D_i)=1$. Similarly, the restriction of $f$ to a curve
$D_i$, $3<i\leq 3+s$, is a bi-double covering of a rational curve
branched at 3 points. Therefore, the geometric genus $g(D_{i})=0$.

(iv) --(v) Since $(f^* (D),f^*(D))_{X_s}=\deg f\cdot (
D,D)_{\widetilde{\mathbb P}^2}=4\cdot ( D,D)_{\widetilde{\mathbb
P}^2}$ for any divisor $D$ on $\widetilde{\mathbb P}^2$, the Claim
follows from the equalities: $(\widetilde
L_j^2)_{\widetilde{\mathbb P}^2}=1-t(L_j)$ for $j=1,\dots,9$ and
$(\widetilde E_i^2)_{\widetilde{\mathbb P}^2}=-1$ for
$i=1,\dots,3+s$. \qed

\begin{picture}(300,260)
\qbezier(40,40)(120,40)(260,40) \qbezier(45,30)(100,140)(155,250)
\qbezier(40,30)(110,110)(180,190)
\qbezier(250,40)(180,120)(117,192)
\qbezier(80,30)(115,135)(150,240)
\qbezier(145,250)(200,140)(255,30)
\qbezier(150,250)(150,140)(150,30)
\qbezier(39,33)(170,120)(218,151)
\qbezier(85,150)(150,107)(260,33) \put(60,88){$L_7$}
\put(152,60){$L_9$} \put(235,80){$L_1$} \put(80,90){$L_5$}
\put(120,77){$L_6$} \put(91,50){$L_8$} \put(205,95){$L_2$}
\put(195,59){$L_3$} \put(181,29){$L_4$} \put(241,30){$p_1$}
\put(50,30){$p_2$} \put(155,240){$p_3$} \put(153,92){$p_4$}
\put(155,154){$p_5$} \put(150,0){$\overline L^{\prime\prime}_2$}
\end{picture}

\begin{center} Fig. 3  \end{center}

Consider  the universal Galois covering $\overline f_s :\overline
X_{s}\to \widetilde{\mathbb P}^2$  and the universal unramified
covering $h_{s,\varphi}:\overline X_s\to X_s$ with respect to
$\varphi :H_1(\mathbb P^2\setminus \overline L_s,\mathbb Z)\to
G=(\mathbb Z/2\mathbb Z)^2$ such that $\overline f_s=f_s\circ
h_{s,\varphi}$. Recall that the covering $\overline f_s$ is
induced by the epimorphism $\psi_{s,\varphi}: H_1(\mathbb P^2
\setminus \overline L_s,\mathbb Z )\to (\mathbb Z /q\mathbb
Z)^{k_{\varphi}}$, where $(\mathbb Z /q\mathbb Z)^{k_{s,\varphi}}$
 and $\psi_{s,\varphi}$ are defined  by
\begin{equation} \label{equi1}
\left\{
\begin{array}{l}
\displaystyle \sum_{j=1}^9 x_j=0, \\
x_{j_1(i)}+x_{j_2(i)}+x_{j_3(i)}=0 , \, \, \, \, \, 3<i\leq 3+s,
\end{array}
\right.
\end{equation}
where for each $i$ the triple $(j_1(i),j_2(i),j_3(i))$ is the set
of indexes of lines $L_j$ such that $p_{i}=L_{j_1(i)}\cap
L_{j_2(i)}\cap L_{j_3(i)}$.

In the case $s=0$ (there are not triple points), we have
$k_{0,\varphi}=8$ and
\begin{equation} \label{s=0}
\deg h_{0,\varphi}=2^6.
\end{equation}

In the case $s=1$, we will assume that $p_{4}=L_{3}\cap L_{6}\cap
L_{9}$ and obtain that $(\mathbb Z /q\mathbb Z)^{k_{1,\varphi}}$
 and $\psi_{1,\varphi}$ are defined  by
\begin{equation} \label{equi2}
\left\{
\begin{array}{l}
\displaystyle \sum_{j=1}^9 x_j=0 \\
x_{3}+x_{6}+x_{9}=0.
\end{array}
\right.
\end{equation}
Therefore, $k_{1,\varphi}=7$ and
\begin{equation} \label{s=1}
\deg h_{1,\varphi}=2^5.
\end{equation}

\begin{picture}(300,260)
\qbezier(46,33)(100,130)(163,224)
\qbezier(250,40)(175,115)(100,190)
\qbezier(80,30)(115,135)(150,240) \qbezier(40,40)(120,40)(260,40)
\qbezier(45,30)(100,140)(155,250)
\qbezier(145,250)(200,140)(255,30)
\qbezier(150,250)(150,140)(150,30)
\qbezier(39,33)(170,120)(218,151)
\qbezier(85,150)(150,107)(260,33) \put(69,105){$L_7$}
\put(78,80){$L_5$} \put(152,165){$L_9$} \put(234,80){$L_1$}
\put(108,100){$L_8$} \put(118,77){$L_6$} \put(197,97){$L_2$}
\put(195,60){$L_3$} \put(181,29){$L_4$} \put(241,30){$p_1$}
\put(50,30){$p_2$} \put(155,240){$p_3$} \put(153,92){$p_4$}
\put(130,165){$p_5$}  \put(150,0){$\overline L^{\prime}_2$}
\end{picture}

\begin{center}  Fig. 4 \end{center}

In the case $s=2$, we will assume that $p_{4}=L_{3}\cap L_{6}\cap
L_{9}$ and $p_5$ (see Fig. 3 and 4) is either the intersection
$L_{2}\cap L_{5}\cap L_{8}$ (a line arrangement $\overline
L_2^{\prime}$) or $L_{2}\cap L_{5}\cap L_{9}$ (a line arrangement
$\overline L_2^{\prime\prime}$).

In the case of a line arrangement $\overline L_2^{\prime}$, we
obtain that $(\mathbb Z /q\mathbb Z)^{k_{2,\varphi}}$  and
$\psi_{2,\varphi}$ are defined  by
\begin{equation} \label{equi31}
\left\{
\begin{array}{l}
\displaystyle \sum_{j=1}^9 x_j=0 \\
x_{3}+x_{6}+x_{9}=0 \\
x_{2}+x_{5}+x_{8}=0
\end{array}
\right.
\end{equation}
and in the case of $\overline L_2^{\prime\prime}$, they are
defined  by
\begin{equation} \label{equi32}
\left\{
\begin{array}{l}
\displaystyle \sum_{j=1}^9 x_j=0 \\
x_{3}+x_{6}+x_{9}=0 \\
x_{2}+x_{5}+x_{9}=0.
\end{array}
\right.
\end{equation}
Therefore in both cases, we have $k_{2,\varphi}=6$ and
\begin{equation} \label{s=2}
\deg h_{2,\varphi}=2^4.
\end{equation}

In the case $s=3$, we will assume that $p_{4}=L_{3}\cap L_{6}\cap
L_{9}$, $p_5=L_{2}\cap L_{5}\cap L_{9}$, and $p_6=L_{2}\cap
L_{6}\cap L_{8}$. We obtain that $(\mathbb Z /q\mathbb
Z)^{k_{3,\varphi}}$  and $\psi_{3,\varphi}$ are defined  by
\begin{equation} \label{equi4}
\left\{
\begin{array}{l}
\displaystyle \sum_{j=1}^9 x_j=0 \\
x_{3}+x_{6}+x_{9}=0 \\
x_{2}+x_{5}+x_{9}=0 \\
x_{2}+x_{6}+x_{8}=0
\end{array}
\right.
\end{equation}
Therefore, $k_{3,\varphi}=5$ and
\begin{equation} \label{s=3}
\deg h_{1,\varphi}=2^3.
\end{equation}

In the case $s=4$, we will assume that $p_{4}=L_{3}\cap L_{6}\cap
L_{9}$, $p_5=L_{2}\cap L_{5}\cap L_{9}$, $p_6=L_{2}\cap L_{6}\cap
L_{8}$, and $p_7=L_{3}\cap L_{5}\cap L_{8}$. The line arrangement
$\overline L_4$ is depicted in Fig. 5.

\begin{picture}(300,270)
\qbezier(50,40)(130,40)(270,40) \qbezier(90,140)(130,140)(230,140)
\qbezier(55,30)(110,140)(165,250)
\qbezier(155,250)(210,140)(265,30)
 \qbezier(160,250)(160,140)(160,30)
\qbezier(49,33)(180,120)(228,151)
\qbezier(95,150)(160,107)(270,33)
\qbezier(165,30)(125,110)(102,156)
\qbezier(155,30)(192,105)(218,156) \put(113,105){$L_1$}
\put(199,105){$L_7$} \put(140,142){$L_4$} \put(116,180){$L_2$}
\put(162,180){$L_9$} \put(192,180){$L_5$} \put(101,77){$L_6$}
\put(205,77){$L_3$} \put(101,29){$L_8$} \put(59,31){$p_{6}$}
\put(249,31){$p_{7}$} \put(147,235){$p_{5}$} \put(161,94){$p_{4}$}
\end{picture}

\begin{center} Fig. 5 \end{center}

We obtain that $(\mathbb Z /q\mathbb Z)^{k_{4,\varphi}}$  and
$\psi_{4,\varphi}$ are defined  by
\begin{equation} \label{equi5}
\left\{
\begin{array}{l}
\displaystyle \sum_{j=1}^9 x_j=0 \\
x_{3}+x_{6}+x_{9}=0 \\
x_{2}+x_{5}+x_{9}=0 \\
x_{2}+x_{6}+x_{8}=0 \\
x_{3}+x_{5}+x_{8}=0
\end{array}
\right.
\end{equation}
It is easy to see that over $\mathbb Z/2\mathbb Z$ the rank of
linear system (\ref{equi5}) is equal to 4. Therefore
$k_{4,\varphi}=5$ and
\begin{equation} \label{s=4}
\deg h_{4,\varphi}=2^3.
\end{equation}

\begin{claim} \label{tors2}
Let $X_{s}$ be a Burniat surface and $\alpha \in \mbox{\rm
Tors}_2(X_s)= \\ \mbox{\rm Tors}_2H^2(X_s,\mathbb Z)$, $\alpha\neq
0$. Then the linear system $|K_{X_s}+\alpha |$ is non-empty and
$D\in |K_{X_s}+\alpha |$ for some $\alpha \in \mbox{\rm
Tors}_2(X_s)$, $\alpha \neq 0$, if and only if $2D= f^*(\overline
D)$ for some $\overline D \in | 3L-\sum_{i=1}^{3+s}E_i|$.
\end{claim}

\proof  The linear system $|K_{X_s}+\alpha |$ is non-empty by
Riemann -- Roch Theorem, since $\dim H^2(X_s, \mathcal
O_{X_s}(K_{X_s}+\alpha))=\dim H^0(X_s, \mathcal
O_{X_s}(\alpha))=0$. Let $D_{\alpha}\in |K_{X_s}+\alpha |$. Then
$2D_{\alpha}\in |2K_{X_s}|$. By Riemann -- Roch Theorem, we have
$$ \dim H^0(X_s,2K_{X_s})=K^2_{X_s}+1=7-s.$$
On the other hand,  by Claim \ref{ample}, $\displaystyle
2K_{X_s}=f^*(3L-\sum_{i=1}^{3+s}E_i)$ and $$ \dim
H^0(\widetilde{P}^2,\mathcal{O}_{\widetilde{P}^2}(3L-\sum_{i=1}^{3+s}E_i))=7-s.$$
Therefore
$$|2K_{X_s}| =f^*(|3L-\sum_{i=1}^{3+s}E_i|)$$
and $D\in |K_{X_s}+\alpha |$ for some $\alpha \in \mbox{\rm
Tors}_2(X_s)$ if and only if $2D= f^*(\overline D)$ for some
$\overline D \in |3L-\sum_{i=1}^{3+s}E_i|$. \qed

\begin{prop} \label{torsb} The $2$-torsion group of a Burniat
surface $X_s$ is isomorphic to $\mbox{\rm Tors}_2(X_s)\simeq
(\mathbb Z/2\mathbb Z)^{6-s}$ if $s\leq 3$ and $\mbox{\rm
Tors}_2(X_4)\simeq (\mathbb Z/2\mathbb Z)^3$.
\end{prop}

\proof It follows from Corollary \ref{tors} that $(\mathbb
Z/2\mathbb Z)^{\deg h_{s,\varphi}}\subset \mbox{\rm Tors}_2(X_s)$.
Note that $\deg h_{s,\varphi}=6-s$  if $s\leq 3$ and $\deg
h_{4,\varphi}=3$. Therefore, by Claim \ref{tors2}, to prove the
Proposition for each case $s=4,3,2,1,0$, it is sufficient to show
that there are exactly $2^{\deg h_{s,\varphi}}-1$ complete
continuous systems of divisors $\overline D$ belonging to $ |
3L-\sum_{i=1}^{3+s}E_i|$ and such that the preimage $f^*(\overline
D)$ of each  $\overline D$  is divisible by two (i.e.,
$f^*(\overline D)=2D$), and each two divisors
$\frac{1}{2}f^*(\overline D_i), \frac{1}{2}f^*(\overline D_j)$
are not linear equivalent if they belong to different systems.

One can check that: \newline in the case $s=4$, the elements
$\overline D \in | 3L-\sum_{i=1}^{7}E_i|$, for which
$f^*(\overline D)$ are divisible by two, are:
\[
\begin{array}{l}
\widetilde L_{3}+\widetilde L_{6}+\widetilde L_{9}+2E_4, \, \, \,
\widetilde L_{2}+\widetilde L_{5}+\widetilde L_{9}+2E_5, \, \, \,
\widetilde L_{2}+\widetilde L_{6}+\widetilde L_{8}+2E_6, \\
\widetilde L_{3}+\widetilde L_{5}+\widetilde L_{8}+2E_7, \, \, \,
\widetilde L_{2}+\widetilde L_{3}+\widetilde L_{7}+E_1, \, \, \,
\,
\, \widetilde L_{5}+\widetilde L_{6}+\widetilde L_{1}+E_2, \\
\widetilde L_{8}+\widetilde L_{9}+\widetilde L_{4}+E_3;
\end{array}
\]
in the case $s=3$, the elements $\overline D \in |
3L-\sum_{i=1}^{6}E_i|$, for which $f^*(\overline D)$ are divisible
by two, are:
\[
\begin{array}{l}
\widetilde L_{3}+\widetilde L_{6}+\widetilde L_{9}+2E_4, \, \, \,
\widetilde L_{2}+\widetilde L_{5}+\widetilde L_{9}+2E_5, \, \, \,
\widetilde L_{2}+\widetilde L_{6}+\widetilde L_{8}+2E_6, \\
\widetilde L_{3}+\widetilde L_{5}+\widetilde L_{8}, \, \, \, \, \,
\, \, \, \, \, \, \, \, \, \, \, \, \, \, \widetilde
L_{2}+\widetilde L_{3}+\widetilde L_{7}+E_1, \, \, \, \,
\, \widetilde L_{5}+\widetilde L_{6}+\widetilde L_{1}+E_2, \\
\widetilde L_{8}+\widetilde L_{9}+\widetilde L_{4}+E_3;
\end{array}
\]
in the case $s=2$ and  $\overline{L}_2=\overline{L}_2^{\prime}$,
the elements $\overline D \in | 3L-\sum_{i=1}^{5}E_i|$, for which
$f^*(\overline D)$ are  divisible by two, are:
\[
\begin{array}{l}
\widetilde L_2+\widetilde L_3+\widetilde L_7+E_1,\, \,\,
\widetilde L_1+\widetilde L_5+\widetilde L_6+E_2,  \, \, \,
\widetilde L_4+\widetilde L_8+\widetilde L_9+E_3, \\
\widetilde L_1+\widetilde L_2+\widetilde L_6+E_1, \, \,  \,
\widetilde L_2+\widetilde L_6+\widetilde L_7+E_2, \, \,  \,
\widetilde L_1+\widetilde L_3+\widetilde L_5+E_1,  \\
\widetilde L_3+\widetilde L_5+\widetilde L_7+E_2,  \, \, \,
\widetilde L_2+\widetilde L_4+\widetilde L_9+E_1,  \,\, \,
\widetilde L_2+\widetilde L_7+\widetilde L_9+E_3, \\
\widetilde L_3+\widetilde L_4+\widetilde L_8+E_1, \, \, \,
\widetilde L_3+\widetilde L_7+\widetilde L_8+E_3, \, \,  \,
\widetilde L_1+\widetilde L_5+\widetilde L_9+E_3, \\
\widetilde L_4+\widetilde L_5+\widetilde L_9+E_2, \, \, \,
\widetilde L_1+\widetilde L_6+\widetilde L_8+E_3, \, \,  \,
\widetilde L_4+\widetilde L_6+\widetilde L_8+E_2;
\end{array}
\]

in the case $s=2$ and
$\overline{L}_2=\overline{L}_2^{\prime\prime}$, the elements
$\overline D \in | 3L-\sum_{i=1}^{5}E_i|$, for which
$f^*(\overline D)$ are  divisible by two,  are:
\[
\begin{array}{l}
\widetilde L_2+\widetilde L_3+\widetilde L_7+E_1, \, \,  \, \, \,
\, \widetilde L_1+\widetilde L_5+\widetilde L_6+E_2,  \, \, \, \,
\, \,
\widetilde L_4+\widetilde L_8+\widetilde L_9+E_3, \\
\widetilde L_1+\widetilde L_2+\widetilde L_6+E_1, \, \,  \, \, \,
\widetilde L_2+\widetilde L_6+\widetilde L_7+E_2, \, \, \, \, \,
\, \,
\widetilde L_1+\widetilde L_3+\widetilde L_5+E_1, \\
\widetilde L_3+\widetilde L_5+\widetilde L_7+E_2, \, \, \, \, \,
\widetilde L_1+\widetilde L_4+\widetilde L_9+E_1, \, \,  \, \, \,
\,
\widetilde L_4+\widetilde L_7+\widetilde L_9+E_2, \\
\widetilde L_2+\widetilde L_5+\widetilde L_9+2E_5, \, \, \,
\widetilde L_3+\widetilde L_6+\widetilde L_9+2E_4, \, \, \,
\widetilde L_2+\widetilde L_6+\widetilde L_8, \\
\widetilde L_1+\widetilde L_7+\widetilde L_9+2E_3, \, \, \,
2\widetilde L_4+\widetilde L_9+E_1+E_2, \, \,  \, \widetilde
L_3+\widetilde L_5+\widetilde L_8.
\end{array}
\]

The reader can test that in the case $s=1$ there are exactly the
$31$ divisors $\overline D \in | 3L-\sum_{i=1}^{4}E_i|$, for which
$f^*(\overline D)$ are divisible by two,  and in the case $s=0$,
there are exactly  $63$ complete continuous systems of divisors
$\overline D$ belonging to $| 3L-\sum_{i=1}^{3}E_i|$, for which
$f^*(\overline D)$ are divisible by two. Note only that in the
case $s=0$, among these systems of divisors, $60$ systems consist
of single divisors and the last $3$ are one-dimensional linear
systems. They are
\[ \widetilde L_1+2\widetilde L_{p_2}+E_2,\, \, \, \widetilde L_4+2\widetilde L_{p_3}+E_3,\, \, \,
\widetilde L_7+2\widetilde L_{p_1}+E_1,
\]
where $\widetilde L_{p_i}=\sigma^{-1}(L_{p_i})$ is the strict
transform of a line belonging to the pencil of lines passing
through the point $p_i$. These three pencils correspond to three
elements, say $\alpha_1, \alpha_2,\alpha_3\in\mbox{\rm
Tors}_2(X_0)$, for which $\dim H^1(X_0,\mathcal
O_{X_0}(\alpha_i))=1$ and these elements "come" from three
irregular intermediate cyclic coverings of the universal Galois
covering $\overline f_0 :\overline X_{0}\to \widetilde{\mathbb
P}^2$ with respect to $\varphi :H_1(\mathbb P^2\setminus \overline
L_0,\mathbb Z)\to G=(\mathbb Z/2\mathbb Z)^2$ (see the end of the
proof of Claim \ref{cl}).

\begin{claim} \label{cl} For $s\geq 1$ the surfaces $\overline X_{s}$
are regular, i.e., the irregularities
$q(\overline X_{s})=0$, and $q(\overline X_{0})=3$.
\end{claim}
\proof The arithmetic genus $p_a$ of a surface is equal to
$p_a=p_g-q+1$. Therefore to calculate $q$, it is sufficient to
calculate $p_a$ and $p_g$.

We have $p_a(X_{s})=1$. Therefore the arithmetic genus
\begin{equation} \label{p-aX-2}
p_a(\overline X_{s})=2^{k_{s,\varphi}-2},
\end{equation}
since $h_s$ is unramified and $\deg h_s=2^{k_{s,\varphi}-2}$.

We calculate the geometric genus $p_g(\overline X_{s})$ for $s\geq
2$ and  the rest two cases are left for the reader, since the
calculation uses the same ideas. As in Section 3, to calculate
$p_g$, it is enough to calculate the geometric genera of
$2^{k_{s,\varphi}}-1$ cyclic coverings corresponding to
$2^{k_{s,\varphi}}-1$ epimorphisms $\psi_m$, $m=1,\dots
,2^{k_{s,\varphi}}-1$, from $G_{u,\varphi}=(\mathbb Z/2\mathbb
Z)^{k_{s,\varphi}}$ to the cyclic group $\mathbb Z/2\mathbb Z$,
where  the group $G_{u,\varphi}$ is isomorphic to the subgroup of
$(\mathbb Z/2\mathbb Z)^9$ given in the coordinates
$(x_1,\dots,x_9)$ by one of linear equations (\ref{equi1}) --
(\ref{equi5}).

It is easy to see that there is a one-to-one correspondence
$\gamma$ between the epimorphisms $\psi_m$ and the elements
$(x_1,\dots x_9)\in G_{u,\varphi}$ such that for $\gamma
(\psi_m)=(x_1,\dots,x_9)$ the cyclic covering corresponding to
$\psi_m$ is given by equation
$$w_m^2=\prod_{x_i=1}l_i.$$
It follows from Claim  \ref{2cov} that the contribution to the
geometric genus of $\overline X_{s}$ can be given only by the
cyclic coverings corresponding to the epimorphisms $\psi_m$ for
which the sum of the coordinates of $\gamma (\psi_m)$ is equal or
grater than $6$, and if it is equal to 6 then the corresponding
branch locus of the cyclic covering has not $4$-fold points.

In the cases $s=3$ or $4$ (see linear equations (\ref{equi4}) and
(\ref{equi5})), we have exactly 7 such coverings given by:
$$
\begin{array}{lll}
z_1^2=l_{1}l_{2}l_{3}l_{5}l_{6}l_{7}  ,\, \, &
z_2^2=l_{1}l_{4}l_{5}l_{6}l_{8}l_{9}  ,\, \, &
z_3^2= l_{2}l_{3}l_{4}l_{7}l_{8}l_{9} , \\
z_4^2=l_{1}l_{2}l_{4}l_{5}l_{7}l_{8}  ,\, \, &
z_5^2=l_{1}l_{3}l_{4}l_{5}l_{7}l_{9} ,\, \, &
z_6^2= l_{1}l_{3}l_{4}l_{6}l_{7}l_{8} ,\\
z_7^2=l_{1}l_{2}l_{4}l_{6}l_{7}l_{9}
\end{array}
$$
Therefore $p_g(\overline X_{3})=p_g(\overline X_{4})=7$ and
$q(\overline X_{3})=q(\overline X_{4})=0$.

In the case $s=2$, when  a line arrangement  $\overline
L_2=\overline L_2^{\prime}$ (see linear equations (\ref{equi31})),
we have exactly 15 such coverings given by:
$$
\begin{array}{lll}
z^2_1=l_{1}l_{2}l_{3}l_{5}l_{6}l_{7}  ,\, \, &
z^2_2=l_{1}l_{2}l_{3}l_{7}l_{8}l_{9}  ,\, \, &
z^2_3= l_{1}l_{4}l_{5}l_{6}l_{8}l_{9} , \\

z^2_4=l_{1}l_{2}l_{4}l_{6}l_{8}l_{9}  ,\, \, &
z^2_5=l_{1}l_{3}l_{4}l_{5}l_{8}l_{9} ,\, \, &
z^2_6= l_{1}l_{2}l_{4}l_{5}l_{6}l_{9} ,\\
z^2_7=l_{1}l_{3}l_{4}l_{5}l_{6}l_{8}  ,\, \, &
z^2_8=l_{1}l_{2}l_{3}l_{6}l_{7}l_{8}  ,\, \, &
z^2_9= l_{1}l_{2}l_{3}l_{5}l_{7}l_{9} , \\
z^2_{10}= l_{2}l_{3}l_{4}l_{5}l_{7}l_{9} ,\, \, &
z^2_{11}=l_{2}l_{3}l_{4}l_{6}l_{7}l_{8}  ,\, \, &
z^2_{12}=l_{2}l_{4}l_{5}l_{6}l_{7}l_{9}  ,\\
z^2_{13}=l_{3}l_{4}l_{5}l_{6}l_{7}l_{8}  ,\, \, &
z^2_{14}=l_{2}l_{4}l_{6}l_{7}l_{8}l_{9} ,\, \, &
z^2_{15}=l_{3}l_{4}l_{5}l_{7}l_{8}l_{9}.
\end{array}
$$
Therefore, $p_g(\overline X_{2}^{\prime})=15$ and $q(\overline
X_{2}^{\prime})=0$.

In the case $s=2$, when  a line arrangement  $\overline
L_2=\overline L_2^{\prime\prime}$ (see linear equations
(\ref{equi32})), we have also exactly 15 such coverings given by:
$$
\begin{array}{lll}
z_1^2=l_{1}l_{2}l_{3}l_{5}l_{6}l_{7}  ,\, \, &
z_2^2=l_{1}l_{4}l_{5}l_{6}l_{8}l_{9}  ,\, \, &
z_3^2= l_{2}l_{3}l_{4}l_{7}l_{8}l_{9} , \\
z_4^2=l_{1}l_{2}l_{4}l_{5}l_{7}l_{8}  ,\, \, &
z_5^2=l_{1}l_{3}l_{4}l_{5}l_{7}l_{9} ,\, \, &
z_6^2= l_{1}l_{3}l_{4}l_{6}l_{7}l_{8} ,\\
z_7^2=l_{1}l_{2}l_{4}l_{6}l_{7}l_{9}
 ,\, \, &
z^2_8=l_{1}l_{3}l_{4}l_{5}l_{7}l_{9}  ,\, \, &
z^2_9= l_{1}l_{2}l_{4}l_{6}l_{8}l_{9} , \\
z^2_{10}= l_{1}l_{3}l_{4}l_{5}l_{8}l_{9} ,\, \, &
z^2_{11}=l_{1}l_{2}l_{3}l_{7}l_{8}l_{9}  ,\, \, &
z^2_{12}=l_{2}l_{3}l_{4}l_{7}l_{8}l_{9}  ,\\
z^2_{13}=l_{2}l_{4}l_{6}l_{7}l_{8}l_{9}  ,\, \, &
z^2_{14}=l_{3}l_{4}l_{5}l_{7}l_{8}l_{9} ,\, \, &
z^2_{15}=l_{1}l_{2}l_{3}l_{4}l_{5}l_{6}l_{7}l_{8}.
\end{array}
$$
The branch locus in the 15-th double covering has degree equal to
$8$ and two 4-fold points. Therefore its geometric genus is equal
to 1. Thus, we have again $p_g(\overline X_{2}^{\prime\prime})=15$
and $q(\overline X_{2}^{\prime\prime})=0$.

The rest two cases are left for the reader, note only that the
non-zero contribution to the irregularity of $\overline X_{0}$ is
given only by the cyclic coverings
$$ z^2_1=l_1l_2l_3l_4, \, \, \, z^2_2=l_1l_7l_8l_9, \, \, \,
z^2_3=l_4l_5l_6l_7.$$ \qed

\begin{cor}
The fundamental group of a Burniat surface $X_0$ is non-abelian
infinite group.
\end{cor}
\proof It follows from Claim \ref{cl}. \qed

It is easy to see that up to a linear transformation of $\mathbb
P^2$ there is the unique Burniat arrangement of lines $\overline
L_4$ (depicted in Fig. 5).

Fix homogeneous coordinates $(z_0:z_1:z_2)$ in $\mathbb P^2$. Put
$\mathcal B_{6-s}=\{ \overline L_s=L_1+\dots+L_9\} $, $s\leq 4$,
the family of the ordered Burniat line arrangements such that
 $L_1$, $L_4$, and $L_7$ are given respectively by $z_0=0$,
$z_1=0$, $z_2=0$, and the point $p_4=(1:1:1)$ (in the case $s=0$
the point $p_4$ is the intersection $L_3\cap L_9$). It is easy to
see that any Burniat arrangement of lines can be transformed to an
arrangement belonging to $\mathcal B_{6-s}$ by a linear
transformation of $\mathbb P^2$. Denote by
$F_s:\mathcal{X}_{6-s}\to \mathcal{B}_{6-s}$ the family of Burniat
surfaces with fibre $F_s^{-1}(\overline L_s)=X_s$ over $\overline
L_s\in \mathcal{B}_{6-s}$, where $X_s$ is the Burniat surface
defined by the line arrangement $\overline L_s$.

If $s=3$ then an arrangement $\overline L_3\in\mathcal B_3$ is
uniquely determined by the point $p_5=L_2\cap L_5\cap L_9 \in
L_9$, since the lines $L_3,L_6,L_9$ are determined by
$p_4=(1:1:1)$, the lines $L_2$ and $L_5$ are determined by $p_5\in
L_9$, and $L_8$ is determined by $p_6=L_2\cap L_6$. Therefore
\begin{equation} \label{B3}
\mathcal B_3\simeq (\mathbb P^1\setminus \{ \mbox{\rm three
points} \})\setminus \mathcal B_2,
\end{equation}
where $\mathcal B_2= \{ \overline L_{4,0}\}$ consists of the
single arrangement $\overline L_{4,0}$ corresponding to the case
$L_3\cap L_5\cap L_8\neq \emptyset$.

Put $\mathcal B_4=\mathcal B_4^{\prime}\cup \mathcal
B_4^{\prime\prime},$ where $\mathcal B_4^{\prime}$ (respectively,
$\mathcal B_4^{\prime\prime}$) consists of the line arrangements
$\overline L_2^{\prime}$ (respectively, $\overline
L_2^{\prime\prime}$).  It is easy to see that an arrangement
$\overline L_2^{\prime}\in \mathcal B_4^{\prime}$ is uniquely
determined by the point $p_5=L_2\cap L_5\cap L_8$. Therefore
\begin{equation} \label{B2'}
\mathcal B_4^{\prime}\simeq \mathbb P^2\setminus \{ \mbox{\rm six
lines} \}.
\end{equation}
Similarly, $\overline L_2^{\prime\prime}\in \mathcal
B_4^{\prime\prime}$ is uniquely determined by the point $p_5\in
L_9$ and the line $L_8$ belonging to the pencil of lines passing
through the point $p_3$. Therefore
\begin{equation} \label{B2''}
\mathcal B_4^{\prime\prime}\simeq (\mathbb P^1\setminus \{
\mbox{\rm three points} \} )^2\setminus (\mathcal B_2 \cup
\mathcal B_3).
\end{equation}

As above, it is clear  that $\overline L_1\in \mathcal B_5$ is
uniquely determined by the lines $L_2$, $L_5$, and $L_8$
belonging, respectively, to the pencils of lines passing through
the points $p_1$, $p_2$, and $p_3$. Therefore
\begin{equation} \label{B1}
\mathcal B_5\simeq (\mathbb P^1\setminus \{ \mbox{\rm three
points} \} )^3\setminus (\mathcal B_2 \cup \mathcal B_3 \cup
\mathcal B_4).
\end{equation}
Similarly,
\begin{equation} \label{B0}
\mathcal B_6\simeq (\mathbb P^1\setminus \{ \mbox{\rm three
points} \} )^4\setminus \overline{\mathcal B}_5,
\end{equation}
where the variety $\overline{\mathcal B}_5$ is the union of
degenerations of the arrangements $\overline L_0$, $\dim
\overline{\mathcal B}_5=3$, $\mathcal B_2 \cup \mathcal B_3 \cup
\mathcal B_4 \cup \mathcal B_5\subset \overline{\mathcal B}_5$.

\begin{claim} \label{noniso} Any two Burniat  surfaces $X_2^{\prime}$ and
$X_2^{\prime\prime}$ are not isomorphic to each other.
\end{claim}

\proof Assume that there is an isomorphism $h:X_2^{\prime} \to
X_2^{\prime\prime}$. Then it induces isomorphisms $h^*: \mbox{\rm
Tors}_2(X_2^{\prime\prime})\to \mbox{\rm Tors}_2(X_2^{\prime})$
and $$h^*:  H^0(X_2^{\prime\prime}, \mathcal
O_{X_2^{\prime\prime}}(K_{X_2^{\prime\prime}}+\alpha)) \to
H^0(X_2^{\prime}, \mathcal
O_{X_2^{\prime}}(K_{X_2^{\prime}}+h^*(\alpha)))$$ for each
$\alpha\in \mbox{\rm Tors}_2(X_2^{\prime\prime})$. Note that Claim
\ref{cl} and Riemann -- Roch Theorem imply $\dim
H^0(X_2^{\prime\prime}, \mathcal
O_{X_2^{\prime\prime}}(K_{X_2^{\prime\prime}}+\alpha))=1$ for
$\alpha\neq 0$. Therefore we should have
$$ K_{h^*(\alpha)}^{\prime}= h^*(K_{\alpha}^{\prime\prime})$$ for
each $K_{\alpha}^{\prime\prime}\in |K_{X_2^{\prime\prime}}+\alpha
| $. On the other hand, among the irreducible components of the
divisors $K_{\alpha}^{\prime\prime}$, there is a rational curve
with self-intersection number equal to $-2$ (the curve $C_9$, see
Claim \ref{k2} and Proposition \ref{torsb}), but among the
irreducible components of the divisors $K_{\alpha}^{\prime}\in
|K_{X_2^{\prime}}+\alpha | $, there is no such a curve.
Contradiction. \qed

\begin{claim} \label{noniso1} For each $s=0,\dots,4$ and for each
$\overline L_{s,0}\in \mathcal{B}_{6-s}$, there are only finitely
many arrangements $\overline L_s\in \mathcal B_{6-s}$ for which
the corresponding Burniat surfaces $X_s=F_s^{-1}(\overline L_s)$
are isomorphic to $X_{s,0}=F_s^{-1}(\overline L_{s,0})$.
\end{claim}

\proof. Put $x=\frac{z_1}{z_0}$ and $y=\frac{z_2}{z_0}$, where
$(z_0:z_1:z_2)$ are the homogeneous coordinates in $\mathbb P^2$
chosen above. Then the lines $L_4$, $L_3$, $L_2$, $L_7$, $L_9$,
$L_8$, $L_5$, and $L_6$ are given, respectively, by equations:
$x=0$, $x=1$, $x=a_1$, $y=0$, $y=1$, $y=b_1$, $x=c_1y$, and
$x=c_2y$ for some $a_1,b_1,c_1,c_2\in \mathbb C\setminus \{
0,1\}$. Consider the injective map $r_s:\mathcal B_{6-s} \to
\mathbb C^{\dim \mathcal B_{6-s}}$ given as follows (if $s=2$ we
will consider two maps: $r_2^{\prime\prime}$ and $r_2^{\prime}$).

In the case $\mathcal B_3$, we have $c_2=1$, $c_1=a_1$, and
$b_1=a_1$, and $a_1$ is a coordinate in $\mathcal B_3$. Put
$r_3(\overline L_3)=a_1\in \mathbb C^1$. Note also that the
arrangement $\mathcal B_2=\{ \overline L_{4,0}\}$ has the
coordinate $a_1=-1$.

In the case $\mathcal B_4^{\prime}$, we have $c_2=1$,
$a_1=b_1c_1$, and $(b_1,c_1)$ are coordinates in $\mathcal
B_4^{\prime}$. Put $r_2^{\prime}(\overline
L_2^{\prime})=(b_1,c_1)\in \mathbb C^2$.

Similarly, in the case $\mathcal B_4^{\prime\prime}$, we have
$c_2=1$, $c_1=a_1$, and $(b_1,c_1)$ are coordinates in $\mathcal
B_4^{\prime\prime}$. Put $r_2^{\prime\prime}(\overline
L_2^{\prime\prime})=(b_1,c_1)\in \mathbb C^2$.

In the case $\mathcal B_5$, we have $c_2=1$, and $(a_1,b_1,c_1)$
are coordinates in $\mathcal B_5$. Put $r_1(\overline
L_1)=(a_1,b_1,c_1)\in \mathbb C^3$.

In the case $\mathcal B_6$, $(a_1,b_1,c_1,c_2)$ are coordinates in
$\mathcal B_6$, and we put $r_0(\overline
L_0)=(a_1,b_1,c_1,c_2)\in \mathbb C^4$.

Obviously, the image $r_s(\mathcal B_{6-s})$ is everywhere dense
open subset of $\mathbb C^{\dim \mathcal B_{6-s}}$.

Let $\mbox{\rm Iso}\, (X_{s,0})$ be the set of arrangements
$\overline L_s$ such that the surfaces $F_s^{-1}(\overline L_s)$
are isomorphic to $X_{s,0}$. Note that $\mbox{\rm Iso}\,
(X_{s,0})$ is a quasi-projec\-tive subvariety of $\mathcal
B_{6-s}$. Indeed, each $X_s$ is a surface with ample canonical
class (possibly, modulo "$-2$"-curves). The imbedding of the
surfaces $X_s$ to $\mathbb P^{10K^2_{X_s}}$ given by $|5K_{X_s}|$
defines a morphism $\mu : \mathcal B_{6-s}\to \mbox{\rm
Hilb}_{P_{X_{s}}}$ to the Hilbert scheme of the surfaces in
$\mathbb P^{10K^2_{X_s}}$ with fixed Hilbert polynomial. The group
$\mbox{\rm PGL}\,(10K^2_{X_s}+1,\mathbb C)$ acts on $\mbox{\rm
Hilb}_{P_{X_{s}}}$ and
$$\mbox{\rm Iso}\, (X_{s,0})=\mu^{-1}(\mu(\mathcal B_{6-s})\cap
\mbox{\rm PGL}(10K^2_{X_s}+1,\mathbb C)(\mu (\overline
L_{s,0})).$$

Therefore, to prove Claim \ref{noniso1}, it is sufficient to show
that the image $r_s(\mbox{\rm Iso}\, (X_{s,0}))$ consists of a
finite set of points.

For this let us consider two isomorphic Burniat surfaces
$X_{s,0}$,  $X_{s,1}$, and let $h: X_{s,0}\to X_{s,1}$ be an
isomorphism. As in the proof of Claim \ref{noniso}, the
isomorphism $h$ induces isomorphisms $h^*: \mbox{\rm
Tors}_2(X_{s,1})\to \mbox{\rm Tors}_2(X_{s,0})$ and
$$h^*:  H^0(X_{s,1}, \mathcal
O_{X_{s,1}}(K_{X_{s,1}}+\alpha)) \to H^0(X_{s,0}, \mathcal
O_{X_{s,0}}(K_{X_{s,0}}+h^*(\alpha)))$$ for each $\alpha\in
\mbox{\rm Tors}_2(X_{s,1})$. Claim \ref{cl} and Riemann -- Roch
Theorem imply $\dim H^0(X_{s,1}, \mathcal
O_{X_{s,1}}(K_{X_{s,1}}+\alpha))=1$ for each $\alpha\neq 0$ if
$s\geq 1$ and for almost all $\alpha\neq 0$ except three
particular values of $\alpha$ if $s=0$. Therefore we should have
$$ K_{0,h^*(\alpha)}= h^*(K_{1,\alpha})$$ for
each $K_{1,\alpha}\in |K_{X_{s,1}}+\alpha | $ (in the case $s=0$
we consider only 60 elements $\alpha$ for which $\dim H^0(X_{s,1},
\mathcal O_{X_{s,1}}(K_{X_{s,1}}+\alpha))=1$). In notation
(\ref{notation}), each divisor $K_{1,\alpha}$ is a linear
combination of the curves $C_{j,1}$, $j=1,\dots,9$, and $D_{i,1}$,
$i=1,\dots, 3+s$. Therefore $h^*(R_1)=R_0$, where
$R_k=\sum_{j=1}^9C_{j,k}+\sum_{i=1}^{3+s}D_{i,k}$ for $k=0,1$, and
the invariants of the curves $D_{i,k}$ and $C_{j,k}$ are
invariants of  surfaces $X_{s,k}$. In particular, the set
$R'(X_{s,k})$, consisting of the components of $R_k$ having
positive genus, is also an invariant. Note that $D_{i,k}\in
R'(X_{s,k})$ for $i=1,2,3$ and $C_{j,k}\in R'(X_{s,k})$ if
$t(L_{j,k})\leq 2$, and, in particular,
$C_{1,k},C_{4,k},C_{7,k}\in R'(X_{s,k})$.

Let $C$ be an elliptic curve. Denote by $B_C$ the subset of
$\mathbb C\setminus \{0,1\} $ consisting of the complex numbers
$c$ such that the curve $C$ can be represented as a two-sheeted
covering $f:C\to \mathbb P^1$ branched at four points
$0,1,c,\infty$. It is well known that for each elliptic curve $C$
the set $B_C$ is finite.

Similarly, for hyperelliptic curve $C$, $g(C)=2$, denote by $B_C$
the subset of $(\mathbb C\setminus \{0,1\})^3 $ consisting of the
triples $(c_1,c_2,c_3)$ of complex numbers, $c_i\neq c_j$ for
$i\neq j$, such that the curve $C$ can be represented as a
two-sheeted covering $f:C\to \mathbb P^1$ branched at six points
$0,1,c_1,c_2,c_3,\infty$. As in the case of elliptic curves, the
set $B_C$ is finite for each curve $C$ of genus two.

Consider the restriction of the covering map $f: X_{s}\to
\widetilde{\mathbb P}^2$ to a curve $C\in R'(X_{s})$. It is a
two-sheeted covering of $\mathbb P^1$. Put
$$ B_{X_{s}}= \bigcup_{C\in R'(X_{s})}B_C.$$
It follows from the above discussion that
\begin{itemize}
\item[(*)] {\it For a Burniat surface $X_s$, the set $ B_{X_{s}}$
is finite, and it is an invariant of $X_s$ up to isomorphism.}
\end{itemize}

Now to complete the proof of Claim \ref{noniso1}, it is sufficient
to notice that:
\newline in the case $\mathcal B_3$, the image
$r_3(\overline L_{3,0})=a_1\in B_{X_{3,0}}$, since $a_1\in
B_{C_{7,0}}$;
\newline in the case $\mathcal B_4^{\prime}$
(and, similarly, in the case $\mathcal B_4^{\prime\prime}$), the
image $r_2(\overline L_{2,0})=(b_1,c_1)$ for some
 $b_1,c_1\in B_{X_{2,0}}$, since $b_1\in B_{C_{4,0}}$
and $c_1\in B_{D_{2,0}}$;
\newline in the case $\mathcal B_5$, the
image $r_1(\overline L_{1,0})=(a_1,b_1,c_1)$ for some
 $a_1,b_1,c_1\in B_{X_{1,0}}$, since $a_1\in
B_{C_{7,0}}$, $b_1\in B_{C_{4,0}}$, and $c_1\in B_{D_{2,0}}$;
\newline in the case $\mathcal B_6$, the
image $r_0(\overline L_{0,0})=(a_1,b_1,c_1,c_2)$ for some
 $b_1$ and $(a_1,c_1,c_2)\in B_{X_{0,0}}$, since  $b_1\in B_{C_{4,0}}$, and $(a_1,c_1,c_2)\in
 B_{D_{9,0}}$.

 \qed

Denote by $\Theta _{X_s}$ the tangent sheaf and by
$\Omega^i_{X_s}$ the sheaf of $i$-differential forms  on $X_s$.

\begin{prop} \label{moduli} For $0\leq s\leq 4$
\begin{itemize}
\item[(i)] $\dim H^0(X_s,\Theta _{X_s})=0$ ;
\item[(ii)]  $\dim H^1(X_s,\Theta
_{X_s})=2s -2+ 3\max (0,2-s)$;
\item[(iii)] $\dim H^2(X_s,\Theta _{X_s})=3\max (0,2-s)$ .
\end{itemize}
\end{prop}
\proof It is known that $\dim H^0(X_s,\Theta _{X_s})=0$ for
surfaces of general type, since the automorphism group of a
surface of general type is a discrete group.

By Riemann -- Roch Theorem, the Euler characteristic $\chi
(\Theta_{X_s})$ of the sheaf $\Theta_{X_s}$ is equal to
$$\chi
(\Theta_{X_s})=\sum_{i=0}^2(-1)^i\dim H^i(X_s,\Theta
_{X_s})=2K^2_{X_s}-10=2s-2,$$ and by Serre duality, $\dim
H^i(X_s,\Theta _{X_s})=\dim H^{2-i}(X_s,\Omega^1_{X_{s}}\otimes
\Omega^2_{X_{s}})$. Therefore to prove Proposition \ref{moduli},
it is sufficient to prove the following Proposition.
\begin{prop} \label{h2} For a Burniat surface $X_s$, $s=0,\dots, 4$, 
$$\dim H^0(X_s,\Omega^1_{X_s}\otimes \Omega^2_{X_s})=3\max
(0,2-s).$$ 
\end{prop}

\proof Put $X=X_s$. Choose a chart $U=\mathbb C^2\subset P^2$ such
that all singular points of the line arrangement $\overline L$ lie
in $\mathbb C^2$, and let $x,y$ be non-homogeneous coordinates in
$U$, $l_i(x,y)=0$ an equation of $L_i\subset \overline L$.

The inclusion of the function field $ \mathbb C(\widetilde{\mathbb
P}^2)\subset\mathbb C(X)$ induced by $f$ is the Galois extension
with Galois group $G=(\mathbb Z/2\mathbb Z)^2$. Let
$\alpha_1=(1,0),\, \alpha_2=(0,1)$, and $\alpha_3=(1,1)$ be the
non-zero elements of $G$. Identifying $\mathbb
C(\widetilde{\mathbb P}^2)$ with the field $K_0=\mathbb C(x,y)$ of
rational functions in $x,y$ and the field $\mathbb C(X_2)$ with
$K=\mathbb C(x,y,w_1,w_2,w_3)$, where $x,y,w_1,w_2,w_3$ satisfy
the equations (\ref{w1w2w3}), without loss of generality, we can
assume that $\alpha (x)=x$, $\alpha (y)=y$ for all $\alpha \in G$
and
$$\begin{array}{l}
\alpha_1(w_1)=w_1,\, \, \alpha_2(w_1)=\alpha_3(w_1)=-w_1, \\
\alpha_2(w_2)=w_2,\, \, \alpha_1(w_2)=\alpha_3(w_2)=-w_2, \\
\alpha_3(w_3)=w_3,\, \, \alpha_1(w_3)=\alpha_2(w_3)=-w_3, \\
\end{array}
$$
Put $K_i=\mathbb C(x,y,w_i)$ and denote by  $g_i:Z_i\to \mathbb
P^2$ the covering induced by the extension $K_0\subset K_i$.

Consider the spaces $M$, $M_0$, $M_1$, $M_2$, $M_3$  of rational
$(1,0)\otimes (2,0)$-forms on $X$, $\mathbb P^2$, $Z_1$, $Z_2$,
$Z_3$, respectively. We have $M_0\subset M_i\subset M$ for
$i=1,2,3$ and
$$\begin{array}{ll}
M_0=\mathbb C(x,y)dx\otimes (dx\wedge dy) \oplus \mathbb
C(x,y)dy\otimes (dx\wedge dy), \\
M=Kdx\otimes (dx\wedge dy) \oplus Kdy\otimes (dx\wedge
dy)\end{array}
$$
and
$$M_i=K_idx\otimes
(dx\wedge dy) \oplus K_idy\otimes (dx\wedge dy)$$ for $i\geq 1$.
Moreover, $G$ acts on $M$ and $M^G=M_0$. Besides, $\omega \in M$
belongs to $M_i$, $i=1,2,3$, if and only if
$\alpha_i(\omega)=\omega$ and $\alpha_j(\omega)=-\omega$ for
$j\neq i$.

The Galois group $G$ acts also on the space
$H^0(X,\Omega^1_{X}\otimes \Omega^2_{X})$ and this space is also
decomposed in the direct sum of eigen-spaces $H_{(i)}$:
$$\displaystyle H^0(X,\Omega^1_{X}\otimes
\Omega^2_{X}) =\displaystyle \bigoplus_{i=0}^3 H_{( i)},$$ where
$\omega \in H_{(i)}$, $i\geq 1$, if and only if
$\alpha_i(\omega)=\omega$ and $\alpha_j(\omega)=-\omega$ for
$j\neq i$, and $\omega \in H_{(0)}$ if and only if
$\alpha_i(\omega)=\omega$ for all $i$. It is easy to see that
\begin{equation}
H_{(i)}=H^0(X,\Omega^1_{X}\otimes \Omega^2_{X})\cap M_i
\end{equation}
for $i=0,1,2,3$.
\begin{lem} \label{local} Let $(V,o)\subset \mathbb C^2\times \mathbb C^1$
be a germ  a surface given in coordinates $(z_1,z_2,w_1)$ by
equation $w_1^2=z_1$.  Consider the action of $\mathbb Z/2\mathbb
Z$ on $H^0(V,\Omega^1_{V}\otimes \Omega^2_{V})$ given by $\alpha
(z_1)=z_1$, $\alpha (z_2)=z_2$, $\alpha (w_1)=-w_1$, where $\alpha
\in \mathbb Z/2\mathbb Z$ is the non-zero element.
\begin{itemize}
\item[(i)] If $\omega \in H^0(V,\Omega^1_{V}\otimes \Omega^2_{V})$
is invariant under the action of $\mathbb Z/2\mathbb Z$, then
$$\omega = (P(z_1,z_2)\frac{dz_1}{z_1}
+Q(z_1,z_2)dz_2)\otimes (dz_1\wedge dz_2);$$ \item[(ii)] if
$\omega \in H^0(V,\Omega^1_{V}\otimes \Omega^2_{V})$ is
anti-invariant under the action of $\mathbb Z/2\mathbb Z$, then
$$\omega = (P(z_1,z_2)dz_1+Q(z_1,z_2)dz_2)\otimes \frac{dz_1\wedge
dz_2}{w_1}, $$
\end{itemize}
where $P(z_1,z_2)$ and $Q(z_1,z_2)$ are some analytic functions in
$z_1$ and $z_2$.
\end{lem}

\proof  Note that $\displaystyle \frac{dz_1\wedge dz_2}{w_1}$ is a
holomorphic nowhere vanishing two-form on $V$ and $\displaystyle
\alpha (\frac{dz_1\wedge dz_2}{w_1})=- \frac{dz_1\wedge
dz_2}{w_1}$. Therefore, a form $\omega\in
H^0(V,\Omega^1_{V}\otimes \Omega^2_{V})$ can be written in the
form
$$\omega =(h_2 dz_2+h_3dw_1)\otimes \frac{dz_1\wedge
dz_2}{w_1},$$ where $h_2,h_3\in H^0(V,\mathcal O_{V})$. Similarly,
the functions $h_i$ can br written in the form `
$h_i=H_i^{\prime}(z_1,z_2) +w_1 H_i^{\prime\prime}(z_1,z_2)$,
where $H_i^{\prime}(z_1,z_2)$ and $H_i^{\prime\prime}(z_1,z_2)$
are some analytic functions in $z_1$ and $z_2$.

It is easy to see that $\omega $ is invariant if and only if
$H_2^{\prime}(z_1,z_2)=H_3^{\prime\prime}(z_1,z_2)=0$, and
 $\omega $ is anti-invariant if and only if
$H_2^{\prime\prime}(z_1,z_2)=H_3^{\prime}(z_1,z_2)=0$. To complete
the proof of Lemma \ref{local}, note that
$w_1dw_1=\frac{1}{2}dz_1$. \qed

\begin{claim} \label{H(o)} The space $ H_{(0)}=0$.
\end{claim}

\proof Let $\omega \in H_{(0)}$, $\omega\neq 0$. It follows from
Lemma \ref{local} that over the chart $U=\mathbb C^2$  it can be
written in the form $\omega= \omega_1\otimes(dx\wedge dy)$, where
$\omega_1\in H^0(U\setminus \mbox{Sing}\,\overline L,
\Omega^1_{U\setminus \mbox{Sing}\,\overline L}(\log \overline L))
$ is $1$-form with logarithmic poles along $\overline L$. Assume
that it has poles along lines $L_{i_1},\dots, L_{i_k}$, $0\leq
k\leq 9$. Then $\omega$ can be written in the form
\begin{equation} \label{omega} \omega=\frac{P(x,y)dx +Q(x,y)dy}{l_{i_1}\dots
l_{i_k}}\otimes (dx\wedge dy), \end{equation} where
$P(x,y),Q(x,y)\in \mathbb C [x,y]$ are polynomials of degree less
or equal $k-4$. Indeed, $l_{i_1}\dots l_{i_k}\omega$ is a regular
form in $\mathbb C^2\setminus \mbox{Sing}\,\overline L$. Therefore
it can be written in the form
$$l_{i_1}\dots l_{i_k}\omega=(P(x,y)dx
+Q(x,y)dy)\otimes (dx\wedge dy),$$ where $P(x,y),Q(x,y)\in \mathbb
C [x,y]$. Next, $\omega$ is regular at the generic point of the
line at infinity $L_{\infty}=\mathbb P^1\setminus C^2$. Let
$(z_0:z_1:z_2)$ be homogeneous coordinates in $\mathbb P^2$ such
that $x=\frac{z_1}{z_0}$ and $y=\frac{z_2}{z_0}$. Then in
coordinates $u=\frac{1}{x},v=\frac{y}{x}$ the form $\omega$ has
the form
$$\omega=((\frac{u^{k-2}\widetilde P(u,v)}{u^{\deg P}\widetilde l_{i_1}\dots \widetilde
l_{i_k}}+\frac{u^{k-2}v\widetilde Q(u,v)}{u^{\deg Q}\widetilde
l_{i_1}\dots \widetilde l_{i_k}})du - \frac{u^{k-1}\widetilde
Q(u,v)}{u^{\deg Q}\widetilde l_{i_1}\dots \widetilde
l_{i_k}}dv)\otimes (\frac{du\wedge dv}{u^3}).$$
Therefore, $\deg
Q\leq k-4$ if $\omega$ is regular at the generic point of the line
at infinity. Similarly, we obtain $\deg P\leq k-4$ if we consider
coordinates $u_1=\frac{1}{y},v_1=\frac{x}{y}$. Therefore, $k\geq
4$.

Let $l_i(x,y)=y+a_ix+b_i$. We have $a_i\neq a_j$ if $i\neq j$,
since all singular points of $\overline L$ lie in $U$, and
\begin{equation} \label{dy} dy=dl_i-a_idx. \end{equation}
Substituting (\ref{dy}) into (\ref{omega}), we obtain
\begin{equation} \label{omega1} \omega=\frac{(P(x,y)-a_iQ(x,y))dx +Q(x,y)dl_i}{l_{i_1}\dots
l_{i_k}}\otimes (dx\wedge dy). \end{equation} Since
$$\frac{(P-a_iQ)dx +Qdl_i}{l_{i_1}\dots
l_{i_k}}\in H^0(U\setminus \mbox{Sing}\,\overline L,
\Omega^1_{U\setminus \mbox{Sing}\,\overline L}(\log \overline
L_2)), $$ the polynomials $P(x,y)-a_{i_j}Q(x,y)$ should be
divisible by $l_{i_j}(x,y)$ for $j=1,\dots, k$. Therefore $k\geq
5$ and the pencil $P(x,y)-aQ(x,y)=0$ of plane curves of degree
$d=k-4$ should have $k$ different fibres containing lines and, by
assumption, each of  these lines is not a fixed component of the
pencil.

Let us show that it is impossible in our case. Indeed, the case
$d=1$ (i.e., $k=5$)  is impossible, since only four lines from
$\overline L$ can lie in the same pencil. The case $d=2$ (i.e.,
$k=6$)  is impossible, since a pencil of conics can have only
three fibers containing lines.

To show that the case $d\geq 3 $ (i.e., $k\geq 7$) is impossible,
note that if the arrangement $l_{i_1}\dots l_{i_k}=0$ has a
$4$-fold point, then the orders of zero of $P(x,y)$ and $Q(x,y)$
at the $4$-fold point $p$ should be at least two. Therefore the
order of zero of  each member of the pencil $P(x,y)-aQ(x,y)=0$ at
$p$ should be also at least two. Indeed, assume that the
arrangement $l_{i_1}\dots l_{i_k}=0$ has such a point $p$. Without
loss of generality, we can assume that $p$ has the coordinates
$(0,0)$. Let $\sigma : \widetilde U\to U$ be the blow up with
center at $p$. In one of the charts, $\sigma$ is given by
equations $x=x_1$ and $y=x_1y_1$. In these new coordinates
$\omega$ has the form
\[
\begin{array}{l}
\omega= \\
\frac{(P(x_1,x_1y_1)+y_1Q(x_1,x_1y_1))dx_1
+x_1Q(x_1,x_1y_1)dy_1} {x_1^4\widetilde l_{i_1}\dots \widetilde
l_{i_k}}\otimes (x_1dx_1\wedge dy_1), \end{array}
\] and therefore the order of
zero of $Q(x,y)$ at $p$ should be at least two, since $\omega$ is
a form with logarithmic poles along the exceptional divisor
$x_1=0$ according to Claim \ref{local}. Similar arguments
(consider the map given by $x=x_2y_2$ and $y=y_2$) show that the
order of zero of $P(x,y)$ at $p$ should be at least two also.

Let us show that the cases $d=3,4,5$ (i.e., $k=7,8,9$) are also
impossible, since in each of these cases, the pencils
$P(x,y)-aQ(x,y)=0$ of degree $d$ should have a fixed line
belonging to the arrangement $l_{i_1}\dots l_{i_k}=0$. Indeed,
each common point of any two lines (belonging to different fibers
of the pencil $P(x,y)-aQ(x,y)=0$) of the arrangement $l_{i_1}\dots
l_{i_7}=0$ is a base point of the pencil. But, it is easy to check
that if we remove any two lines from $\overline L$ (the case
$d=3$, i.e., $k=7$), then we obtain a new arrangement consisting
of seven lines such that there is a component of the new
arrangement passing through four its singular points (counting
with multiplicities). Therefore this line should be a fixed
component of the pencil of degree 3. Similarly,  if we remove any
line from $\overline L$ (the case $d=4 $, .e., $k=8$),  then we
obtain a new arrangement consisting of eight lines such that there
is a component of the new arrangement passing through a $4$-fold
point and three other singular points of the new arrangement. In
the case $d=5$ (i.e., $k=9$), also there is a line (for example,
$L_1$) passing through two $4$-fold points and two other singular
points of the line arrangement $\overline L$. \qed
\begin{claim} \label{H(i)} Let $X$ be a Burniat surface $X_s$. Then
$$\dim H_{(1)}=\dim H_{(2)}=\dim H_{(3)}=\max (0,2-s).$$
\end{claim}

\proof Consider the space $H_{(1)}$ (the cases of the spaces
$H_{(2)}$ and $H_{(3)}$ are similar). Let $\omega \in H_{(1)}$,
$\omega\neq 0$. Then, since $\omega \in M_1$, we have
$$\omega = (R_1(x,y)dx+R_2(x,y)dy)\otimes \frac{dx\wedge
dy}{w_1}, $$ where $R_i(x,y)$ are rational functions. Note that
the form $\frac{dx\wedge dy}{w_1}$ has not poles and zeros on
$Y_1\setminus \mbox{Sing}\, Y_1$, since $w_1^2=l_{1}\dots l_{6}$
and $\deg l_{1}\dots l_{6}=6$ (see Section 3). Therefore, since
$\omega$ is a regular form over the generic point of $L_{\infty}$,
as in the proof of Claim \ref{H(o)}, applying Lemma \ref{local},
one can easily show that $\omega$ can be written in the form
\begin{equation} \label{omega11} \omega=\frac{P(x,y)dx +Q(x,y)dy}{l_{7}l_{8}l_{9}}
\otimes \frac{dx\wedge dy}{w_1},
\end{equation}
where $P(x,y),Q(x,y)\in \mathbb C[x,y]$ are polynomials of degree
$\leq 2$, and, moreover, the form
$$\omega_1=\frac{P(x,y)dx +Q(x,y)dy}{l_{7}l_{8}l_{9}}\in H^0(U\setminus \mbox{Sing}\, D_{1},
\Omega^1_{U\setminus Sing\, D_{1}}(\log D_{1})),$$ where
$D_{1}=L_7+L_8+L_9$. Without loss of generality, we can assume
that $p_3=(0,0)$, $l_1(x,y)=x-ay$, $l_7(x,y)=y$, $l_8(x,y)=x-y$,
and $l_9(x,y)=x$, where $a\neq 0,\, 1$.

To resolve singularities of $Y_s$ we should blow up 4-fold and
triple points of $\overline L_s$. If the form $\omega\in H_{(1)}$,
then it should be regular over the blown up curves $E_i$.

The following Lemmas \ref{local2} -- \ref{local4} give necessary
and sufficient conditions on the form (\ref{omega11}) to be
regular over the curve $E_3$.
\begin{lem} \label{local2} Let $(V,o)\subset (\mathbb C^2\times \mathbb C^2,o)$
be a germ of a normal surface given in coordinates
$(z_1,z_2,w_1,w_2)$ by equations $w_1^2=x-ay$ and $w_2^2=xy(x-y)$,
where $a\neq 0,\, 1$, and let
$$\omega=\frac{P(x,y)dx +Q(x,y)dy}{xy(x-y)}\otimes
\frac{dx\wedge dy}{w_1} \in H^0(\overline V, \Omega^1_{\overline
V}\otimes \Omega^2_{\overline V}),$$ where $P(x,y),\, Q(x,y)\in
\mathbb C[x,y]$, $\deg P(x,y)=\deg Q(x,y)=2$, and $\nu: \overline
V\to V$ is the minimal resolution of the singular point $o$ of
$V$. Then
$$
\displaystyle \omega=(c\frac{ydx-xdy}{xy(x-y)}+\frac{P_2(x,y))dx
+Q_2(x,y)dy}{xy(x-y)})\otimes \frac{dx\wedge dy}{w_1},
$$
where $c$ is a constant and $P_2(x,y)$, $Q_2(x,y)$ are homogeneous
polynomials of degree $2$.
\end{lem}

\proof Let $Z$ be the image $g(V)$ of the map
$g((x,y,w_1,w_2))=(x,y)$ and $\sigma :\widetilde Z\to Z$ the blow
up with center at $g(o)$, $E=\sigma^{-1}(g(o))$. By Lemma \ref{lem
1.3}, the map $f:\overline V\to \widetilde Z$ induced by $g$ is an
analytic covering and it can be factorized into the composition
$f=\overline f_1\circ h_1$, where $\overline f_1: \overline V_1\to
\widetilde Z$ is a $\mathbb Z/2\mathbb Z$-covering and $\overline
V_1$ is bi-meromorphic to the surface given by $w_1^2=x-ay$.

The morphism $\sigma$ is given by $x=u,\, y=uv$ in some local
coordinates in $\widetilde Z$. Then the surface $\overline V_1$ is
given locally by $w_1^2=u(1-av)$. Therefore $\overline f_1$ is
branched along $E$ and it is easy to see that $h_1$ is not
branched at the generic point of $\overline f_1^{-1}(E)$. Thus
$h_1$ is a local isomorphism at the generic point of $\overline
f^{-1}(E)$.

The form $\omega$ is a meromorphic form on $\overline V_1$ and, by
assumption, $h^*(\omega)$ is holomorphic. Therefore, $\omega$ is
holomorphic at the generic point of $\overline f_1^{-1}(E)$.
Moreover, by Lemma \ref{local}, it can have at most logarithmic
poles along the curves given by equations $y=0$, $x-y=0$, and
$x-ay=0$. We have
\[ \begin{array}{l} \displaystyle
\omega=\frac{(P(u,uv)+vQ(u,uv))du +uQ(u,uv)dv}{u^3v(1-v)}\otimes
\frac{udu\wedge
dv}{w_1}= \\
\displaystyle \hspace{1cm}(\frac{(P(u,uv)+vQ(u,uv))du}
{u^2v(1-v)}+\frac{ Q(u,uv)dv}{uv(1-v)})\otimes (dw_1\wedge dv).
\end{array}
\]
Therefore $Q(0,0)=P(0,0)=0$ and  $P(u,uv)+vQ(u,uv)$ should be
divisible by $u^2$.

Put $P=a_1x+a_2y+P_2(x,y)$ and $Q=b_1x+b_2y+Q_2(x,y)$, where $P_2$
and $Q_2$ are homogeneous polynomials of degree two. We have
$$P(u,uv)+vQ(u,uv)=a_1u+a_2uv+b_1uv+b_2uv^2+P_2(u,uv)+Q_2(u,uv),$$
where $P_2(u,uv)+Q_2(u,uv)$ is divisible by $u^2$. Therefore
$a_1=b_2=0$ and $a_2=-b_1$. We have
$$
\displaystyle \omega=(a_2\frac{ydx-xdy}{xy(x-y)}+\frac{P_2(x,y))dx
+Q_2(x,y)dy}{xy(x-y)})\otimes \frac{dx\wedge dy}{w_1}.
$$
\qed

\begin{lem} \label{local4} Let $x,y$ be coordinates in $U=\mathbb
C^2$, $o=(0,0)$, $D\subset U$ a divisor given by $xy(x-y)=0$, and
$$\omega_1=\frac{P(x,y)dx +Q(x,y)dy}{xy(x-y)}\in H^0(U\setminus \{ 0\},
\Omega^1_{U\setminus \{ o\}}(\log D)),$$  where $P(x,y)$ and
$Q(x,y)$  are some homogeneous polynomials of degree two. Then
$\omega_1$ is a linear combination of $\frac{dx}{x}$,
$\frac{dy}{y}$, and $\frac{d(x-y)}{x-y}$.
\end{lem}

\proof  Since $\omega_1 \in H^0(U\setminus \{ (0,0)\},
\Omega^1_{U\setminus \{ (0,0)\} }(\log D)$, then $P=yP_1(x,y)$ and
$Q=xQ_1(x,y)$. Let $P_1(x,y)=a_1x+a_2y$ and $Q_1(x,y)=b_1x+b_2y$.
Put $l=x-y$. We have $dl=dx-dy$. Therefore
$$\begin{array}{l}
\displaystyle \omega_1=\frac{yP_1(x,y)dx +xQ_1(x,y)dy}{xy(x-y)}=
\\
\displaystyle \frac{((x-l)P_1(x,x-l)+xQ_1(x,x-l))dx
-xQ_1(x,x-l)dl}{x(x-l)l}= \\
\displaystyle \frac{((a_1+a_2+b_1+b_2)x^2+l(\dots))dx
-xQ_1(x,y)dl}{x(x-l)l}\end{array}
$$
and, consequently, we should have
\[
 a_1+a_2+b_1+b_2=0,
\]
i.e.,
$$\begin{array}{l}
\displaystyle \omega_1=\frac{(a_1xy+a_2y^2)dx
+(b_1x^2-(a_1+a_2+b_1)xy)dy}{xy(x-y)}= \\
\displaystyle \hspace{1cm}
-a_2\frac{dx}{x}+b_1\frac{dy}{y}+(a_1+a_2)\frac{d(x-y)}{x-y}.
\end{array}
$$  \qed

It follows from Lemmas \ref{local2} -- \ref{local4} that if
$\omega\in H_{(1)}$, then $\omega$ has the form
\begin{equation} \label{form1}
\displaystyle \omega=
(c\frac{ydx-xdy}{y(x-y)(x-ay)}+c_1\frac{dx}{x}+c_2\frac{dy}{y}+c_3\frac{d(x-y)}{x-y})\otimes
 \frac{dx\wedge dy}{w_1}.
\end{equation}

\begin{lem} \label{local5} A form
$$
\displaystyle \omega_1=
c_1\frac{dx}{x}+c_2\frac{dy}{y}+c_3\frac{d(x-y)}{x-y}
$$  is regular at the generic point of the line at infinity
$L_{\infty}=\mathbb P^2\setminus \mathbb C^2$
if and only if $c_1+c_2+c_3=0$.
\end{lem}

\proof  Let $x=\frac{1}{u}$ and $y=\frac{v}{u}$. We have
$$\begin{array}{l}
\displaystyle \omega_1=
c_1\frac{dx}{x}+c_2\frac{dy}{y}+c_3\frac{d(x-y)}{x-y}= \\
\displaystyle \hspace{1cm}(c_1+c_2+c_3)\frac{du}{u}+
\frac{c_2(1-v)-c_3v}{v(1-v)}dv .
\end{array}
$$  \qed

\begin{lem} \label{local6}
Led $D\subset \mathbb P^2$ be the projective closure of the curve
$D_0\subset \mathbb C^2$ given by $xy(x-y)=0$ and $0=(0,0)\in
\mathbb C^2$ the origin. Then the form
$$
\displaystyle \frac{ydx-xdy}{xy(x-y)}\in H^0(\mathbb P^2\setminus
\{ 0\}, \Omega^1_{\mathbb P^2 \setminus \{ o\}}(\log D)).
$$
\end{lem}
\proof Straightforward. \qed

It follows from Lemmas \ref{local5} and \ref{local6} that if
$\omega$ having the form (\ref{form1}) belongs to $H_{(1)}$ then
\begin{equation} \label{equation1}
c_1+c_2+c_3=0.
\end{equation}

Let $p_2$ has coordinates $(b,0)$. The following Lemma  gives a
necessary and sufficient condition on the form (\ref{form1}) to be
regular over the curve $E_2$.

\begin{lem} \label{local7} Let  $\overline V$ be the desingularization of a germ
of a surface $(V,o)\subset (\mathbb C^2\times \mathbb C^2,0)$
given in coordinates $(x,y,w_1,w_3)$ by equations
$w_1^2=(x-a_1y)(x-a_2y)(x-a_3y)$ and $w_2^3=y$, and let
$$\omega=
 (c\frac{ydx-(x+b)dy}{(x+b)y(x-y+b)}+c_2\frac{dy}{y})\otimes
 \frac{dx\wedge dy}{w_1}
 \in
H^0(\overline V, \Omega^1_{\overline V}\otimes \Omega^2_{\overline
V}),$$ where $a_1$, $a_2$, $a_3$, and $b$ are constants, $a_i\neq
a_j$ for $i\neq j$, $a_i\neq 0$ for $i=1,2,3$, $b\neq 0$. Then
$c-bc_1=0$.
\end{lem}

\proof  Consider the covering $g:V\to g(V)=Z\subset \mathbb C^2$
given by $g((x,y,w_1,w_3))=(x,y)$ and  let $\sigma
:\widetilde{Z}\to Z$ be the blow up with center at $g(o)=(0,0)$,
$E=\sigma^{-1}(g(o))$. By Lemma \ref{lem 1.3}, the map
$f:\overline V\to \widetilde {Z}$ induced by $g$ is a regular
covering and it can be factorized into the composition
$f=\overline f_1\circ h_1$, where $\overline f_1: \overline V_1\to
\widetilde Z$ is a $\mathbb Z/2\mathbb Z$-covering and $\overline
V_1$ is bimeromorphic to the surface given by
$w_1^2=(x-a_1y)(x-a_2y)(x-a_3y)$.

Let $\sigma$ be given by $x=uv,\, y=v$ in some local coordinates
in $\widetilde Z$. Then the surface $\overline V_1$ is given
locally by $\widetilde w_1^2=v(u-a_1)(u-a_2)(u-a_3)$, where
$\widetilde w_1=\frac{w_1}{v}$ and the surface $\overline V$ is
the normalization of a surface given locally by $\widetilde
w_1^2=v(u-a_1)(u-a_2)(u-a_3)$ and $\widetilde w_3^2=v$. Therefore
$\overline f_1$ is  branched along the exceptional curve $E$ given
by $v=0$ and $h_1$ is not branched at the generic point of
$\overline f_1^{-1}(E)$. Thus, $h_1$ is a local isomorphism at the
generic point of $f^{-1}(E)$.

The form
$$\begin{array}{l}
\displaystyle \omega=
(c\frac{ydx-(x+b)dy}{(x+b)y(x-y+b)}+c_2\frac{dy}{y})\otimes
 \frac{dx\wedge dy}{w_1}= \\
\displaystyle \hspace{1cm}\frac{cv^2du+(-cb+c_2b^2
+c_2v(\dots))dv}{v(uv+b)(uv-v+b)} \otimes \frac{du\wedge
dv}{\widetilde w_1}.
\end{array}
$$
Since $\omega\in H^0(\overline V, \Omega^1_{\overline V}\otimes
\Omega^2_{\overline V})$ and $b\neq 0$, it follows from Lemma
\ref{local} that $c-c_2b$ should be equal to $0$. \qed

It follows from Lemma \ref{local7} that if $\omega$ having the
form (\ref{form1}) belongs to $H_{(1)}$ then
\begin{equation} \label{equation2}
c-bc_2=0,
\end{equation}
where $p_2=(0,b)$.

One can check that the point $p_1$ does not give a restriction on
the form (\ref{form1}).

Consider restrictions on the form (\ref{form1})are given by triple
points.  Let one of the lines $L_8$ or $L_9$ (say $L_9$) passes
through a triple point $p_{3+i}$ of $\overline L_s$. Since the
point $p_{3+i}\in L_9$, it has coordinates $(0,b_1)$, $b_1\neq 0$

\begin{lem} \label{local3} Let $\overline V$ be the desingularization of a germ
of a surface $(V,o)\subset (\mathbb C^2\times \mathbb C^2,o)$ be a
given in coordinates $(x,y,w_1,w_2)$ by equations
$w_1^2=(x-a_1y)(x-a_2y)$ and $w_2^2=x(x-a_1y)$, and let
$$\omega=
 (c\frac{(y+b_1)dx-xdy}{x(y+b_1)(x-y-b_1)}+c_1\frac{dx}{x})\otimes
 \frac{dx\wedge dy}{w_1}
$$ belongs to $H^0(\overline V, \Omega^1_{\overline V}\otimes
\Omega^2_{\overline V})$, where $a_1$, $a_2$, and $b_1$ are
constants, $a_1\neq a_2$, $a_i\neq 0$, $b_1\neq 0$. Then
$c-b_1c_1=0$.
\end{lem}

\proof  As in the proof of Lemma \ref{local2}, consider the map
$g:V\to g(V)=Z\subset \mathbb C^2$ given by
$g((x,y,w_1,w_2))=(x,y)$ and  let $\sigma :\widetilde{Z}^2\to Z$
be the blow up with center at $g(o)=(0,0)$, $E=\sigma^{-1}(g(o))$.
By Lemma \ref{lem 1.3}, the map $f:\overline V\to \widetilde {Z}$
induced by $g$ is a regular covering and it can be factorized into
the composition $f=\overline f_1\circ h_1$, where $\overline f_1:
\overline V_1\to \widetilde Z$ is a $\mathbb Z/2\mathbb
Z$-covering and $\overline V_1$ is bimeromorphic to the surface
given by $w_1^2=(x-a_1y)(x-a_2y)$.

Let $\sigma$ be given by $x=u,\, y=uv$ in some local coordinates
in $\widetilde Z$. Then the surface $\overline V_1$ is given
locally by $\widetilde w_1^2=(1-a_1v)(1-a_2v)$, where $\widetilde
w_1=\frac{w_1}{u}$ and the surface $\overline V$ is the
normalization of a surface given locally by $\widetilde
w_1^2=(1-a_1v)(1-a_2v)$ and $\widetilde w_2^2=(1-a_1v)$, where
$\widetilde w_2=\frac{w_2}{u}$. Therefore $\overline f_1$ is not
branched along the exceptional curve $E$ given by $u=0$ and $h_1$
is not branched at the generic point of $\overline f_1^{-1}(E)$.
Thus $h_1$ is a local isomorphism at the generic point of
$f^{-1}(E)$.

We have
$$\begin{array}{l}
\displaystyle
\omega=(c\frac{(y+b_1)dx-xdy}{x(y+b_1)(x-y-b_1)}+c_1\frac{dx}{x})\otimes
 \frac{dx\wedge dy}{w_1}
= \\
\displaystyle \hspace{1cm}\frac{cb_1-c_1b_1^2 +u(\dots )du
+u(\dots) dv}{u(uv+b_1)(u-uv-b_1)} \otimes \frac{du\wedge
dv}{\widetilde w_1}.
\end{array}
$$
Since $\omega\in H^0(\overline V, \Omega^1_{\overline V}\otimes
\Omega^2_{\overline V})$, $b_1\neq 0$, the number $c-c_1b_1$
should be equal to $0$. \qed

It follows from Lemma \ref{local3} that if $\omega$ having the
form (\ref{form1}) belongs to $H_{(1)}$, then
\begin{equation} \label{equation3}
c-c_1b_1=0,
\end{equation}
where $p_{3+i}=(0,b_1)\in L_9$.

If the arrangement $\overline L$ has two triple points with
coordinates $(0,b_1)$ and $(0,b_2)$ lying in $L_9$, then the
equations $c-c_1b_1=0$ and $c-c_1b_2=0$ are linear independent.
Similarly, it is easy to see that if a triple point $p_{s+3}\in
L_8$, then it also give some linear equation of the form
\begin{equation} \label{equation4}
f(c,c_3)=0. \end{equation}

As a consequence, we obtain that the space $H_{(1)}$ consists of
the forms
\[
\displaystyle \omega=
(c\frac{ydx-xdy}{xy(x-y)}+c_1\frac{dx}{x}+c_2\frac{dy}{y}+c_3\frac{d(x-y)}{x-y})\otimes
 \frac{dx\wedge dy}{w_1}
\]
satisfying $2+s$ linear equations (\ref{equation1}) --
(\ref{equation4}). Note that these equations are linear
independent. Therefore $\dim H_{(1)} =\max (0,2-s)$. \qed \qed

Denote by $\mathcal M_{6-s}$, $0\leq s\leq 4$, the union of
irreducible components of the moduli scheme  of the surfaces of
Burniat type $s$, and denote by  $\widetilde{\mathcal B}_{6-s}$
the image of $\mathcal B_{6-s}$ in $\mathcal M_{6-s}$.

\begin{cor} \label{modBur} The variety $\widetilde{\mathcal B}_{6-s}$ is everywhere
dense in $\mathcal M_{6-s}$ if $s\leq 2$ and
\begin{itemize}
\item[(i)] the space $\mathcal M_2$ is non-singular at $X_{4,0}=\widetilde{\mathcal B}_2$,
$\dim \mathcal M_2=6$;
\item[(ii)] the space $\mathcal M_3$ is non-singular at any point $X_3\in\widetilde{\mathcal B}_3$,
$\dim \mathcal M_3=4$, and  $\widetilde{\mathcal B}_3$ is a
rational curve;
\item[(iii)] the space $\mathcal M_4=M_4^{\prime}\cup
M_4^{\prime\prime}$ consists of two irreducible rational surfaces
$\mathcal M^{\prime}_4$ and $\mathcal M^{\prime\prime}_4$,
$\mathcal M_4$ is non-singular at each point $X_2\in
\widetilde{\mathcal B}_4$;
\item[(iv)]
the space $\mathcal M_5$ is unirational $3$-fold nonsingular at
any point $X_1\in\widetilde{\mathcal B}_5$;
\item[(v)]
the space $\mathcal M_6$ is unirational $4$-fold nonsingular at
any point $X_0\in\widetilde{\mathcal B}_6$.
\end{itemize}
\end{cor}

\proof It follows from (\ref{B3}) -- (\ref{B0}), Claims
\ref{noniso}, \ref{noniso1}, and Proposition \ref{moduli}.

\begin{prop} \label{peters} In notation of Section 4.1,
let the Campedelli covering $f_C:X_C\to \widetilde{\mathbb P}^2$
be branched along the Campedelli line arrangement depicted in Fig.
6. Then the Burniat surface $X_4$ is isomorphic to the Campedelli
surface $X_C$.
\end{prop}

\begin{picture}(300,270)

\qbezier(50,40)(130,40)(270,40) \qbezier(45,140)(130,140)(250,140)
\qbezier(55,30)(110,140)(165,250)
\qbezier(155,250)(210,140)(265,30)

\qbezier(160,250)(160,140)(160,30)
\qbezier(49,33)(180,120)(228,151)
\qbezier(95,150)(160,107)(270,33)

\put(55,144){$L_{(0,1,0)}$} \put(50,85){$L_{(1,1,0)}$}
\put(162,160){$L_{(0,1,1)}$} \put(240,84){$L_{(1,1,1)}$}
\put(99,57){$L_{(1,0,0)}$} \put(192,60){$L_{(1,0,1)}$}
\put(101,29){$L_{(0,0,1)}$} \put(59,31){$p_{1}$}
\put(249,31){$p_{3}$} \put(147,235){$p_{2}$}
\put(117,145){$p_{4}$} \put(209,147){$p_{5}$}
\put(161,94){$p_{6}$}
\end{picture}

\begin{center} Fig. 6 \end{center}

\proof First of all, note that $f_C$ can be decomposed into the
composition $f_C=f_{C,1}\circ h_{C,1}$, where $\overline f_{C,1}:
X_{C,1}\to \widetilde{\mathbb P}^2$ is the desingularization of
the double covering $g_{C,1}:Y_{C,1}\to \mathbb P^2$ given in
non-homogeneous coordinates by equation
$$w_1^2=l_{(1,0,0)}l_{(1,1,0)}l_{(1,0,1)}l_{(1,1,1)}$$
(here $l_{\alpha}(x,y)=0$ is an equation of the line
$L_{\alpha}$). To resolve the singularities of $Y_{C,1}$, we
should blow up the points $p_1\dots,p_6$. Let $\sigma:
\widetilde{\mathbb P}^2 \to \mathbb P^2$ the composition of these
blowups, $E_i=\sigma^{-1}(p_i)$, and
$\widetilde{L}_{\alpha}=\sigma^{-1}({L}_{\alpha})$.

It is easy to see that the curves $E_i$ do not belong to the
branch locus of $f_{C,1}$ and for each $(0,a_2,a_3)$ the strict
transform $f_{C,1}^{-1}(\widetilde{L}_{(0,a_2,a_3)})$ of
$\widetilde{L}_{(0,a_2,a_3)}$ is the disjoint union of two
rational curves $L_{(0,a_2,a_3)}^{\prime}$ and
$L_{(0,a_2,a_3)}^{\prime\prime}$, since the rational curve
$\widetilde{L}_{(0,a_2,a_3)}$ does not meet the branch locus of
$f_{C,1}$. Therefore the $(\mathbb Z/2\mathbb Z)^2$-covering
$h_{C,1}: X_C\to X_{C,1}$ is branched over the union of the curves
$E_i$, $i=1,\dots ,6$, and the curves
$f_{C,1}^{-1}(\widetilde{L}_{(0,a_2,a_3)})$, $(a_2,a_3)\in
(\mathbb Z/2\mathbb Z)^2$. We have
$(L_{(0,a_2,a_3)}^{\prime},L_{(0,a_2,a_3)}^{\prime})_{X_{C,1}}=
(L_{(0,a_2,a_3)}^{\prime\prime},L_{(0,a_2,a_3)}^{\prime\prime})_{X_{C,1}}=-1$,
since the intersection number
$$(\widetilde{L}_{(0,a_2,a_3)},\widetilde{L}_{(0,a_2,a_3)})_{\widetilde{\mathbb
P}^2}=-1$$ and $\deg f_{C,1}=2$. The mutual arrangement of the
curves $L_{(0,a_2,a_3)}^{\prime}$ and
$L_{(0,a_2,a_3)}^{\prime\prime}$, $(a_2,a_3)\in (\mathbb
Z/2\mathbb Z)^2$, is depicted in Fig. 7.

\begin{picture}(50,120)
\put(100,50){\special{em:moveto}}
\put(140,100){\special{em:lineto}}
\put(210,100){\special{em:lineto}}
\put(250,50){\special{em:lineto}} \put(210,0){\special{em:lineto}}
\put(140,0){\special{em:lineto}} \put(100,50){\special{em:lineto}}
\put(165,107){$L_{(0,1,0)}^{\prime}$}
\put(165,-13){$L_{(0,1,0)}^{\prime\prime}$}
\put(230,15){$L_{(0,1,1)}^{\prime}$}
\put(82,73){$L_{(0,1,1)}^{\prime\prime}$}
\put(91,15){$L_{(0,0,1)}^{\prime}$}
\put(234,73){$L_{(0,0,1)}^{\prime\prime}$}
\end{picture}
\vspace{0.5cm}

\begin{center} Fig. 7 \end{center}

Similarly, for each $(1,a_2,a_3)$ the
 intersection number
$$(\widetilde{L}_{(1,a_2,a_3)},\widetilde{L}_{(1,a_2,a_3)})_{\widetilde{\mathbb
P}^2}=-2$$ and therefore, the strict transform
$D_{(1,a_2,a_3)}=f_{C,1}^{-1}(\widetilde{L}_{(1,a_2,a_3)})$ has
the self-intersection number
$$(D_{(1,a_2,a_3)},D_{(1,a_2,a_3)})_{X_{C,1}}=
(\frac{1}{2}f_{C,1}^*(\widetilde{L}_{(1,a_2,a_3)}),
\frac{1}{2}f_{C,1}^*(\widetilde{L}_{(1,a_2,a_3)}))_{X_{C,1}}=-1$$
and
$$(D_{(1,a_2,a_3)},L_{(0,b_2,b_3)}^{\prime})_{X_{C,1}}=
(D_{(1,a_2,a_3)},L_{(0,b_2,b_3)}^{\prime\prime})_{X_{C,1}}=0$$ for
all $(a_2,a_3)$ and $(b_2,b_3)$. Note that each $D_{(1,a_2,a_3)}$
is also a rational curve.

It is not hard to see that $X_{C,1}$ is a rational surface and if
$\tau: X_{C,1}\to \widetilde{X}_{C,1}$ is the blowdown of the
curves $L_{(0,1,0)}^{\prime}$, $L_{(0,0,1)}^{\prime}$,
$L_{(0,1,1)}^{\prime}$, and four curves $D_{(1,a_2,a_3)}$, then
$\widetilde{X}_{C,1}$ is isomorphic to the projective plane
$\mathbb P^2$, the image
$$\tau (\sum_{i=1}^6
E_i+L_{(0,1,0)}^{\prime\prime}+L_{(0,0,1)}^{\prime\prime}+L_{(0,1,1)}^{\prime\prime})$$
is the Burniat line arrangement $\overline L_4$, and the covering
$h_{C,1}$ coincides with the Burniat covering $f:X_4\to
\widetilde{\mathbb P}^2$. \qed

\begin{cor} \label{cambur} The moduli space $\mathcal M_2$ coincides with
the moduli space $\mathcal C$ of the Campedelli surfaces.
\end{cor}
{\bf 4.3} {\it A surface $X$ of general type with $p_g=0$,
$K_X^2=6$ and} $(\mathbb Z/3\mathbb Z)^3\subset \mbox{Tors}(X)$.
Let $\overline L=L_1+\dots +L_6$ be an arrangement in $\mathbb
P^2$ of six lines having $3$ triple points $p_1,\, p_2,\, p_3$ not
lying in the same line. The arrangement $\overline L$ is depicted
in Fig. 8.

\begin{picture}(300,270)

\qbezier(50,40)(130,40)(270,40) \qbezier(55,30)(110,140)(165,250)
\qbezier(155,250)(210,140)(265,30)
\qbezier(160,260)(160,150)(160,30)
\qbezier(50,32)(130,99)(210,166) \qbezier(95,150)(160,107)(270,33)
\put(47,45){$p_2$} \put(265,45){$p_1$} \put(165,240){$p_3$}
\put(116,180){$L_1$} \put(162,180){$L_6$} \put(192,180){$L_2$}
\put(101,87){$L_5$} \put(205,77){$L_4$} \put(101,29){$L_3$}
\end{picture}

\begin{center} Fig. 8 \end{center}

Consider a covering $g: Y\to \mathbb P^2$ associated with the
epimorphism $\varphi :H_1(\mathbb P^2\setminus \overline L,
\mathbb Z)\to G=(\mathbb Z/3\mathbb Z)^2$ given  by
\[
\begin{array}{l}
\varphi (\lambda_1)=\varphi (\lambda_2)=\varphi (\lambda_3)=(1,0),
\\ \varphi (\lambda_4)=(2,1),\, \, \varphi (\lambda_5)=(1,1),\,
\, \varphi (\lambda_6)=(0,1).
\end{array}
\]
The surface $Y$ has $3$ singular points lying over the triple
points $p_i$. By Lemma \ref{lem 1.3} to resolve them, it is
sufficient to blow up the points $p_i$ and consider the induced
Galois covering $f: X \to \widetilde{\mathbb P}^2$, where $\sigma
: \widetilde{\mathbb P}^2\to {\mathbb P}^2$ is the composition of
blowups with centers at the points $p_i$. Denote by
$E_i=\sigma^{-1}(p_i)$ the exceptional curve lying over $p_i$.

\begin{prop} \label{newex} The constructed above surface $X$ is a
surface of general type with $K_X^2=6$, $p_g=0$, and $(\mathbb
Z/3\mathbb Z)^3\subset \mbox{Tors}(X)$.
\end{prop}

\proof By Claim \ref{ample}, we have $3K_X=|f^*(3L-\sum E_i)|$,
where $L=\sigma^*(\mathbb P^1)$ is the total transform of a line
$\mathbb P^1\subset \mathbb P^2$. Therefore $X$ is a surface of
general type with ample canonical class. Applying (\ref{K2ex}) and
(\ref{eXex}), it is easy to see that $K_X^2=6$ and $e(X)=6$.
Therefore, by Noether's formula, $p_a=1-q+p_g=1$. As above, to
calculate $p_g$, it is enough to calculate the geometric genera of
$4$ cyclic coverings corresponding to $4$ epimorphisms from the
group $G=(\mathbb Z/3\mathbb Z)^2$ to the cyclic group $\mathbb
Z/3\mathbb Z$. These coverings are given respectively in
non-homogeneous coordinates by the following equations:
\[
\begin{array}{l}
w_1^3=l_1l_2l_3l_4^2l_5; \, \,  w_2^3=l_1l_2l_3l_5^2l_6; \\
w_3^3=l_1l_2l_3l_4l_6^2;\, \,  w_4^3=l_4l_5l_6.
\end{array}
\]
Applying Claim \ref{3cov}, one can easily check that the geometric
genus of each of these coverings is equal to zero. Thus, $X$ has
the geometric genus $p_g=0$.

To see that $(\mathbb Z/3\mathbb Z)^3\subset \mbox{Tors}(X)$,
consider the universal covering $g_{u(3)}: \widetilde Y_{u(3)}\to
\widetilde{\mathbb P}^2$ corresponding to the epimorphism
$$\overline
\varphi :H_1(\mathbb P^2\setminus \overline L, \mathbb Z)\to
H_1(\mathbb P^2\setminus \overline L, \mathbb Z/3\mathbb Z )\simeq
(\mathbb Z/3\mathbb Z)^5$$ given by
\[
\begin{array}{ll}
\varphi (\lambda_1)=(1,0,0,0,0), \, \, & \varphi
(\lambda_2)=(1,0,1,0,0), \\ \varphi (\lambda_3)=(1,0,0,1,0),\, \,
& \varphi (\lambda_4)=(2,1,0,0,1), \\ \varphi
(\lambda_5)=(1,1,0,0,0),\, \, & \varphi (\lambda_6)=(0,1,2,2,2).
\end{array}
\]
It is easy to see that the Galois covering $h_{u,\varphi}:X_u\to
X$ induced by the projection $\psi :(\mathbb Z/3\mathbb Z)^5\to
G=(\mathbb Z/3\mathbb Z)^2$ to the first two coordinates is
unramified.  Therefore, by Corollary \ref{tors}, $(\mathbb
Z/3\mathbb Z)^3\subset \mbox{Tors}{X}$

\begin{claim} \label{cl33} The surface $X_{u}$  has irregularity
$q(X_{u})=3$.
\end{claim}
\proof As in the proof of Claim \ref{cl}, to calculate $q$, it is
sufficient to calculate $p_a$ and $p_g$.

We have $p_a(X)=1$. Therefore, the arithmetic genus
$p_a(X_{u})=3^3$, since $h_{u,\varphi}$ is unramified and $\deg
h_{u,\varphi}=3^3$.

To calculate $p_g$, it is sufficient to calculate, by Claim
\ref{3cov}, the geometric genera of $\frac{3^{5}-1}{2}=121$ cyclic
coverings corresponding to $\frac{3^{5}-1}{2}$ epimorphisms
$\psi_m$, $m=1,\dots ,121$, of $G_{u,\varphi}=(\mathbb Z/3\mathbb
Z)^{5}$ to the cyclic group $\mathbb Z/3\mathbb Z$. The
calculation is left to the reader, note only that the contribution
to the irregularity of $X_{u}$ is given only by the cyclic
coverings
$$ z^3_1=l_1l_2l_6, \, \, \, z^3_2=l_1l_3l_5, \, \, \,
z^3_3=l_2l_3l_4.$$ \qed

\begin{cor}
The fundamental group of the surface $X$, constructed above, is a
non-abelian infinite group.
\end{cor}
\proof It follows from Claim \ref{cl33}. \qed

{\bf 4.4.} {\it The Godeaux surface}. Let $\overline
L=L_1+L_2+L_3+L_4$ be an arrangement in $\mathbb P^2$ of four
lines in general position. Consider the following coverings: the
universal covering $g_{u(5)}: Y_{u(5)}\to \mathbb P^2$
corresponding to the epimorphism
$$\overline
\varphi :H_1(\mathbb P^2\setminus \overline L, \mathbb Z)\to
H_1(\mathbb P^2\setminus \overline L, \mathbb Z/5\mathbb Z )\simeq
(\mathbb Z/5\mathbb Z)^3,$$ a covering $g: Y\to \mathbb P^2$
associated with the epimorphism $\varphi :H_1(\mathbb P^2\setminus
L, \mathbb Z)\to (\mathbb Z/5\mathbb Z)^2$ given in some chosen
coordinates in $G=(\mathbb Z/5\mathbb Z)^2$ by
$$\varphi (\lambda_1)=(1,0),\, \, \varphi (\lambda_2)=(0,1),\, \,
\varphi (\lambda_3)=(1,2),\, \, \varphi (\lambda_4)=(3,2),$$ and
the covering $h:Y_{u(5)}\to Y$ corresponding to an epimorphism
$\psi :(\mathbb Z/5\mathbb Z)^3\to G=(\mathbb Z/5\mathbb Z)^2$
such that $\varphi=\psi\circ\overline \varphi$. By Lemma \ref{lem
1.3}, the surface $Y$ is nonsingular and by Proposition \ref{uni},
the covering $h$ is unramified.
\begin{prop} \label{Godo} The constructed above surface $Y$ is a
surface of general type with $K_Y^2=1$, $p_g=0$, and $\mbox{\rm
Tors}\, (Y)=\mathbb Z/5\mathbb Z$.
\end{prop}

\proof By Claim \ref{ample}, we have $5K_Y=|f^*(L)|$, where $L$ is
a line in $\mathbb P^2$. Therefore $Y$ is a surface of general
type with ample canonical class. Applying (\ref{K2ex}) and
(\ref{eX}), it is easy to see that $K_Y^2=1$ and $e(Y)=11$.
Therefore, by Noether's formula, $p_a=1-q+p_g=1$. To calculate
$p_g$, it is enough to calculate the geometric genera of  $6$
cyclic coverings corresponding to $6$ cyclic subgroups of $G$ and
given respectively in non-homogeneous coordinates by the following
equations:
\[
\begin{array}{l} w_1^5=l_1l_3l_4^3;\, \,
w_2^5=l_2l_3^2l_4^2;\, \, w_3^5=l_1l_2l_3^3; \\
w_4^5=l_1^2l_2l_3^4l_4^3;\, \, w_5^5=l_1l_2^2l_4^2;\, \,
w_6^5=l_1l_2^3l_3^2l_4^4.
\end{array}
\]
Applying calculation made in section 3, one can easily check that
the geometric genus of each of these coverings is equal to zero.
Thus, $Y$ has the geometric genus $p_g=0$.

Since the Galois covering $h$ is unramified, we have the following
inclusion: $\mathbb Z/5\mathbb Z \subset \mbox{Tors}\, (Y)$. To
show that $\mbox{Tors}\, (Y)=\mathbb Z/5\mathbb Z $, it is
sufficient to show that $Y_{u(5)}$ is simply connected. Moreover,
it is easy to see that $Y_{u(5)}$ is isomorphic to a smooth
surface in $\mathbb P^3$. Indeed, let us choose homogeneous
coordinates $(x_0:x_1:x_2)$ in $\mathbb P^2$ such that $x_i=0$ is
an equation of $L_{i+2}$. Let $\sum a_ix_i=0$ be an equation of
$L_1$. Without loss of generality we can assume that the covering
$g_{u(5)}$ is associated to the epimorphism $\overline \varphi
:H_1(\mathbb P^2\setminus L, \mathbb Z)\to (\mathbb Z/5\mathbb
Z)^3$ given by
$$\overline \varphi (\lambda_1)=(0,0,1),\, \,
\overline \varphi (\lambda_2)=(4,4,4),\, \, \overline \varphi
(\lambda_3)=(1,0,0),\, \, \overline \varphi (\lambda_4)=(0,1,0).$$
In this case $Y_{u(5)}$  is given by
$$ z_3^5=z_1;\, \, z_4^5=z_2;\, \, z_5^5=a_0+a_1z_1+a_2z_2$$ in non-homogeneous coordinates
$(z_1,z_2,\dots,z_5)$, where $z_1=\frac{x_1}{x_0}$ and
$z_2=\frac{x_2}{x_0}$, and therefore $Y_{u(5)}$ is isomorphic to
the projective closure of the surface in $\mathbb C^3$ given by
$z_5^5=a_0+a_1z_3^5+a_2z_4^5$ (cf. \cite{God}).

\ifx\undefined\bysame
\newcommand{\bysame}{\leavevmode\hbox to3em{\hrulefill}\,}
\fi

\end{document}